\documentclass[11pt, reqno]{amsart}
\usepackage{physics}
\usepackage{verbatim}
\usepackage{mcode}
\usepackage{listings}
\usepackage{graphicx}
\usepackage{float}
\usepackage[utf8]{inputenc}
\usepackage[margin=1in]{geometry}

\theoremstyle{plain}
\newtheorem{theorem}{Theorem}[section]

\newtheorem{cor}{Corollary}[section]

\newcommand{\la}{\lambda}

\newcommand{\eps}{\epsilon}

\newenvironment{proof1}[1]{\begin{trivlist} \item[] {\em Proof of #1:}}{\hfill $\Box$
                      \end{trivlist}}
                      
\DeclareRobustCommand{\rchi}{{\mathpalette\irchi\relax}}
\newcommand{\irchi}[2]{\raisebox{\depth}{$#1\chi$}} 

\graphicspath {{figures/}}

\makeindex

\title[Limiting eigenfunctions of a Spectral Flow]{Limiting Eigenfunctions of Sturm-Liouville operators Subject to a Spectral Flow}

\author[T. Beck]{Thomas Beck}
\email{tdbeck@email.unc.edu}
\address{Department of Mathematics, University of North Carolina at Chapel Hill \\ CB\#3250
  Phillips Hall \\ Chapel Hill, NC 27599}
  
  \author[I. Bors]{Isabel Bors}
\email{isabelpbors@gmail.com}
\address{Department of Mathematics, University of North Carolina at Chapel Hill \\ CB\#3250
  Phillips Hall \\ Chapel Hill, NC 27599}
  
    \author[G. Conte]{Grace Conte}
\email{gconte23@live.unc.edu}
\address{Department of Mathematics, University of North Carolina at Chapel Hill \\ CB\#3250
  Phillips Hall \\ Chapel Hill, NC 27599}
  
  \author[G. Cox]{Graham Cox}
  \email{gcox@mun.ca}
\address{Department of Mathematics and Statistics, Memorial University of Newfoundland, St. John's, NL A1C 5S7, Canada}

\author[J.L. Marzuola]{Jeremy L. Marzuola}
\email{marzuola@math.unc.edu}
\address{Department of Mathematics, University of North Carolina at Chapel Hill \\ CB\#3250
  Phillips Hall \\ Chapel Hill, NC 27599}

\begin{document}

\begin{abstract}
We examine the spectrum of a family of Sturm--Liouville operators with regularly spaced delta function potentials parametrized by increasing strength. The limiting behavior of the eigenvalues under this spectral flow was described in \cite{nodal}, where it was used to study the nodal deficiency of Laplacian eigenfunctions. Here we consider the eigenfunctions of these operators. In particular, we give explicit formulas for the limiting eigenfunctions, and also characterize the eigenfunctions and eigenvalues for all values for the spectral flow parameter (not just in the limit). We also develop spectrally accurate numerical tools for comparison and visualization.   
\end{abstract}

\maketitle

\section{Introduction}
It is well known that the $n$-th eigenfunction of a regular Sturm--Liouville problem, with separable boundary conditions on a finite interval, has precisely $n-1$ interior zeros. Put differently, it has exactly $n$ nodal domains, where the nodal domains of $u$ are defined to be the connected components of the set $\{x : u(x) \neq 0\}$.

In higher dimensions one has the Courant nodal domain theorem, which says the $n$th eigenfunction of the Laplacian, or more generally the Schr\"odinger operator $-\Delta + V$, has \emph{at most} $n$ nodal domains. Unlike the one-dimensional case, this is generally a strict inequality: While Courant's theorem says that $\nu(n) \leq n$ for all $n$, where $\nu(n)$ denotes the number of nodal domains for the $n$-th eigenfunction, for any dimension greater than one it was shown by Pleijel \cite{P56} that the equality $\nu(n) = n$ can only hold for finitely many values of $n$.

An eigenfunction for which $\nu(n) = n$ is said to be \emph{Courant sharp}. For some relative simple domains, such as squares, balls, equilateral triangles and tori, it is possible to completely determine the Courant sharp eigenfunctions. See, for instance, \cite{berard2020courant,helffer2016nodal,berard2016courant,bonnaillie2015nodal,lena2015courant,helffer2010spectral,helffer2009nodal} among many others.  In general this is a very challenging problem, and there is much that is not known.

The extent to which an eigenfunction is not Courant sharp is measured by its \emph{nodal deficiency},
\begin{equation}
    \delta(n) := n - \nu(n).
\end{equation}
The first explicit formula for the nodal deficiency that we are aware of was given in \cite{BKS12}, in terms of the Morse index of an energy function defined on the space of equipartitions. A second formula, in terms of Dirichlet-to-Neumann maps on the nodal set, was obtained in \cite{CJM_nodal} using infinite-dimensional symplectic methods.

The starting point of the current work is \cite{nodal}, in which the nodal deficiency formula of \cite{CJM_nodal} was reinterpreted (and reproved) using a spectral flow approach. The benefit of this approach is that it identifies an explicit mechanism by which low energy eigenfunctions can contribute to the nodal deficiency.

The idea is as follows:
\begin{enumerate}
    \item Fix an eigenfunction $u_n$, with eigenvalue $\lambda_n$, and define the nodal set $\mathcal{Z}= \overline{\{x : u_n(x) = 0\}}$;
    \item Define the family of operators $L(\sigma) = -\Delta + \sigma \delta_{\mathcal{Z}}$ for $\sigma \geq 0$;
    \item Write the eigenvalues of $L(\sigma)$ in analytic branches $\lambda_m(\sigma)$ such that $\lambda_m(0)$ are the eigenvalues of $L(0) = -\Delta$.
\end{enumerate}
It is not hard to show that the eigenvalue curves $\lambda_m(\sigma)$ are non-decreasing. In particular, $\lambda_n(\sigma)$ remains constant, whereas $\lambda_1(\sigma), \ldots, \lambda_{n-1}(\sigma)$ are strictly increasing. Of these first $n$ eigenvalue curves, precisely $\nu(n)$ of them will converge to $\lambda_n$ as $\sigma \to \infty$, while the remaining $n - \nu(n)$ will converge to values larger than $\lambda_n$. Therefore, \emph{the nodal deficiency of $u_n$ is equal to the number of eigenvalue curves that pass through $\lambda_n$ as $\sigma$ ranges from $0$ to $\infty$.}

In \cite{nodal} this spectral flow was analyzed for a Schr\"odinger operator with separable potential on a rectangular domain, assuming the eigenfunction $u_n$ of interest is separable. These assumptions reduce the problem to a one-dimensional one, where the analysis is easier. The main simplification is due to the fact that in one dimension the eigenvalues must be simple, and so the eigenvalue curves cannot intersect as $\sigma$ changes. (In particular, this implies that in one dimension the nodal deficiency is always zero, reproducing the classic result of Sturm.)

While the analysis of \cite{nodal} completely describes the limiting behavior of the eigenvalues in the separable case, the corresponding eigenfunctions were not studied, as they were not needed to determine the nodal deficiency. The current paper addresses this issue.

This is motivated by the desire to study small perturbations of separable eigenfunctions, for an eigenvalue of multiplicity at least two. For instance, the methods of \cite{nodal} apply to the Laplacian eigenfunction $\sin (2\pi x) \sin (\pi y)$ on the unit square, but not to the eigenfunction $\sin(2\pi x) \sin (\pi y) + \epsilon \sin (\pi x) \sin (2\pi y)$. In the first case the nodal set is a union of horizontal and vertical lines, so the operator $L(\sigma)$ is separable, while in the latter case the nodal set is more complicated, and a similar ODE reduction is not possible. It should be possible to study the spectral flow for such a non-separable eigenfunction by viewing it as a small perturbation of the separable case, but this will require knowledge of the limiting eigenfunctions in the separable case.

Moreover, understanding the limiting eigenfunctions should lead to additional insight into the spectral flow in general, and the mechanism by which an eigenfunction $u_m$ with eigenvalue $\lambda_m < \lambda_n$ does (or does not) contribute to the nodal deficiency of $u_n$.

In order to state our results, we first recall the one-dimensional spectral flow from \cite{nodal}. Consider the
Sturm--Liouville eigenvalue problem
\begin{equation}\label{Vdiffeq}
    -u''(x) + V(x) u = \lambda u(x), \qquad u(0) = u(1) = 0.
\end{equation}
Given the $n$-th eigenfunction $u_n$ of \eqref{Vdiffeq}, the nodal set $\{x \in (0,1) : u_n(x) = 0\}$ consists of $n-1$ points, and hence can be written as $\{ x_k \}_{k=1}^{n-1}$.  For a fixed $n > 1$, the spectral flow is defined by the 
eigenvalue problem
\begin{equation}\label{Vdiffsum}
    -u''(x) + V(x) u(x) + \sigma u(x) \sum_{k=1}^{n-1} \delta(x - x_k)  = \lambda (\sigma) u(x), \qquad u(0) = u(1) = 0,
\end{equation}
parameterized by $\sigma \in [0,\infty)$. Here the potential $V$ has been shifted by adding delta functions of strength $\sigma$ at the zeros of $u_n$.

An equivalent formulation, without reference to delta functions, is that $u(x)$ solves the differential equation \eqref{Vdiffeq} on each open subinterval $(x_{k-1}, x_k)$ (letting $x_0 = 0$ and $x_n =1$ for convenience) with the boundary conditions $u(0) = u(1) = 0$ and
\begin{equation}
    u(x_k^+) = u(x_k^-), \qquad u'(x_k^+) - u'(x_k^-) = \sigma u(x_k) \label{Vdeltabc}
\end{equation}
for $1 \leq k \leq n-1$, where $\pm$ superscripts denote right- and left-hand limits, respectively.
These boundary conditions ensure that the eigenfunctions are continuous at $x_k$, but have a jump in their derivatives whenever $\sigma>0$ and $u(x_k) \neq 0$.

The one-dimensional analysis in \cite{nodal} shows that the eigenvalues $\lambda_m(\sigma)$ of \eqref{Vdiffsum} depend analytically on $\sigma$, with
\[
    \lim_{\sigma \to \infty} \lambda_m(\sigma) = \lambda_n, \qquad 1 \leq m \leq n,
\]
and
\[
    \lim_{\sigma \to \infty} \lambda_m(\sigma) > \lambda_n, \qquad m > n.
\]
This result holds for any continuous potential $V$, and does not rely on explicit formulas for the eigenvalues, which are only available for special choices of $V$. However, it gives no information about the eigenvalues for intermediate values of $\sigma$, or about the eigenfunctions for any values of $\sigma$ other than $0$. In this paper we completely solve this problem for the case $V(x) = 0$.

To that end, we consider the eigenvalue problem
\begin{equation}\label{diffeq}
-u''(x) = \lambda u(x), \ \  \
u(0) = u(1) = 0.
\end{equation}
The eigenfunctions are given by $u_n(x)$ = $C_n\sin(n\pi x)$, where $C_n$ is a constant, with corresponding positive eigenvalues $\lambda_n$ = $n^2\pi^2$, for each $n \in \mathbb{N}$.

We fix $n\geq2$, and impose the boundary conditions \eqref{Vdeltabc} at the zeros of the eigenfunction $u_n(x)$ = $\sin(n\pi x)$. The zeros of $u_n$ occur at $x_k = \tfrac{k}{n}$, where 0 $\leq$ $k$ $\leq$ $n$, $k$ $\in$ $\mathbb{Z}$.  
At the interior zeros, the boundary conditions from \eqref{Vdeltabc} therefore become 
\begin{equation}\label{deltabc}
    u\big(\tfrac{k}{n}^+\big) = u\big(\tfrac{k}{n}^-\big), \qquad
    u'\big(\tfrac{k}{n}^+\big) - u'\big(\tfrac{k}{n}^-\big) = \sigma u\big(\tfrac{k}{n}\big)
\end{equation}
for $1 \leq k \leq n-1$.

For each $\sigma \geq0$, we want to examine the behavior of the eigenvalues $\la_m(\sigma)$ and eigenfunctions $u_m(x;\sigma)$ of the equation and boundary conditions given in \eqref{diffeq} and \eqref{deltabc}. In this notation $\la_m(0) = m^2\pi^2$, $u_m(x;0) = \sin(m\pi x)$, and the eigenvalues $\la_m(\sigma)$ are simple for any finite $\sigma$. Since $u_n(x) = \sin(n\pi x)$ vanishes at the interior nodes $x_k$, the $n$-th eigenvalue is constant, $\la_n(\sigma) = n^2\pi^2$ for all $\sigma$, and we can also take $u_n(x;\sigma)$ to be independent of $\sigma$. However, the same is not true for $\la_m(\sigma)$ and $u_m(x;\sigma)$ for $1 \leq m \leq n-1$. Our main theorem provides a description of these eigenvalues and eigenfunctions for all $\sigma\geq0$.

\begin{theorem}\label{thm:main}
For each $\sigma \geq0$ and $1 \leq m \leq n-1$, the eigenvalues $\la_m(\sigma) = \gamma_m(\sigma)^2$ satisfy the implicit equation
\begin{equation}
\sigma = 2 \gamma_m(\sigma) \frac{\cos\left(\tfrac{m\pi}{n}\right) - \cos\big(\tfrac{\gamma_m(\sigma)}{n}\big)}{\sin\big(\tfrac{\gamma_m(\sigma)}{n}\big)}.
\end{equation}
Set $I_k = [x_{k-1},x_k]$, with $x_k = \tfrac{k}{n}$. Then, for each $\sigma\geq0$ and $1 \leq m \leq n-1$, up to an overall normalization factor, the eigenfunctions $u_m(x;\sigma)$ are given by
\begin{equation}\label{umexplicit}
    u_m(x;\sigma) = \sin\big(\gamma_m(\sigma) x\big) + \sum_{j=1}^{k-1}A_{j,m}(\sigma) \sin\big(\gamma_m(\sigma)(x- x_j)\big) \quad \text{for } x \in I_k,
\end{equation}
where the coefficients $A_{j,m}(\sigma)$ are determined in terms of the eigenvalue $\la_m(\sigma)$ via
\begin{align}
\begin{split}\label{Akmexplicit}
A_{1,m}(\sigma) & = 2\cos\left(\tfrac{m\pi}{n}\right) - 2\cos\big(\tfrac{\gamma_m(\sigma)}{n}\big) ,\\ 
A_{k,m}(\sigma) & = \frac{\sin \left(  \tfrac{k m \pi}{n} \right) }{ \sin \left(  \tfrac{m \pi}{n} \right) } A_{1,m}(\sigma) 
\quad\text{for } 2 \leq k \leq n-1 .
\end{split}
\end{align} 
\end{theorem}

To obtain the limiting eigenfunctions, which we denote by
\begin{equation}\label{def:uminfty}
u_{m}(x;\infty) = \lim_{\sigma\to\infty}u_{m}(x;\sigma),
\end{equation}
one can use the fact that $\gamma_m(\sigma) \to n \pi$ for $1 \leq m \leq n$ to obtain
\[
    u_{m}(x;\infty) =  \left(1 + \sum_{j=1}^{k-1}(-1)^jA_{j,m}\right)\sin(n \pi x)  \quad \text{for } x \in I_k
\]
from \eqref{umexplicit}, where the coefficients $A_{j,m} = \lim_{\sigma\to\infty} A_{j,m}(\sigma)$ are given by \eqref{Akmexplicit}. An alternate approach is given below in Corollary \ref{cor:infty}. Along the way we obtain more information about the eigenfunctions, which leads directly to an explicit formula for $u_m(x;\infty)$, see \eqref{conjsum2} and \eqref{Bkmexplicit}.

As $\sigma$ increases, the derivatives of $u_m(x;\sigma)$ remain bounded, and so to ensure that the interior condition in \eqref{deltabc} continues to hold, the values $u_m(x_k;\sigma)$ must converge to $0$ as $\sigma$ converges to infinity. Our first corollary of Theorem \ref{thm:main} is that these values converge to $0$ at the same rate for each node $x_k$.
\begin{cor} \label{cor:nodes}
Up to an overall normalization factor, for each $\sigma \geq0$ and $1 \leq m \leq n-1$, the values of the eigenfunctions $u_m(x;\sigma)$ at the interior nodes $x_k$ are given by
\begin{align*}
u_m(x_k;\sigma) = \frac{\sin\big(\tfrac{\gamma_m(\sigma)}{n}\big)}{\sin\left(\tfrac{m\pi}{n}\right)}\sin\left(\tfrac{km\pi}{n}\right) = \frac{\sin\big(\tfrac{\gamma_m(\sigma)}{n}\big)}{\sin\left(\tfrac{m\pi}{n}\right)}u_m(x_k;0).
\end{align*}
\end{cor}

We next define the quantity
\begin{equation}
F_{k,m}(\sigma) = \int_{I_k}u_{m}(x;\sigma)\sin(n \pi x) \,d x,
\end{equation}
which provides a measure of the weight of $u_{m}(x;\sigma)$ in each interval $I_k$. Using Theorem \ref{thm:main}, we can carefully control its dependence on $\sigma$.

\begin{cor} \label{cor:integral}
For each $\sigma\geq0$ and $1 \leq m \leq n-1$, the integrals $F_{k,m}(\sigma)$ are equal to
\begin{align*}
F_{k,m}(\sigma) = n\pi(-1)^{k+1}\frac{u_{m}(x_k;\sigma) + u_{m}(x_{k-1};\sigma)}{n^2\pi^2 - \la_m(\sigma) }.
\end{align*}
In particular, choosing a normalization of $u_{m}(x;\sigma)$ so that $F_{1,m}(\sigma) \equiv 1$, the integrals $F_{k,m}(\sigma)$ are then independent of $\sigma$ for $1\leq k \leq n$.
\end{cor}

Finally, we examine the behavior of the eigenfunctions as the parameter $\sigma$ approaches infinity.  As shown in \cite{nodal}, the eigenvalues $\la_m(\sigma)$ for $1 \leq m \leq n$ all converge to $\la_n = n^2\pi^2$ as $\sigma$ tends to infinity.  (Note that this is consistent with our implicit expression for the eigenvalues $\la_{m}(\sigma)$ from Theorem \ref{thm:main}.) From Corollary \ref{cor:nodes}, this ensures that $u_m(x_k;\sigma)$ converges to zero as $\sigma$ tends to infinity. This means that $u_m(x;\infty)$ (defined in \eqref{def:uminfty}) is proportional to $\sin(n\pi x)$ on each interval $I_k$, so it can be represented by a vector with $n$ entries, where the $k$th entry of the vector is the coefficient of $\sin(n \pi x)$ on $I_k$. Our final corollary of Theorem \ref{thm:main} gives an explicit expression for these vectors.

\begin{cor} \label{cor:infty}
For each $m$, $1 \leq m \leq n$, up to an overall normalizing factor, the limiting eigenfunctions $u_{m}(x;\infty)$ are given by
\begin{align} \label{conjsum2}
    u_{m}(x;\infty) =  B_{k,m} \sin(n\pi x) \quad \text{for } x \in I_k,
\end{align}
where 
\begin{equation}\label{Bkmexplicit}
    B_{k,m} = (-1)^{k+1} \sin{\frac{(2k-1)m\pi}{2n}}.
\end{equation}
\end{cor}

\begin{figure}[ht]
  \centering
    \includegraphics[width=.25\textwidth]{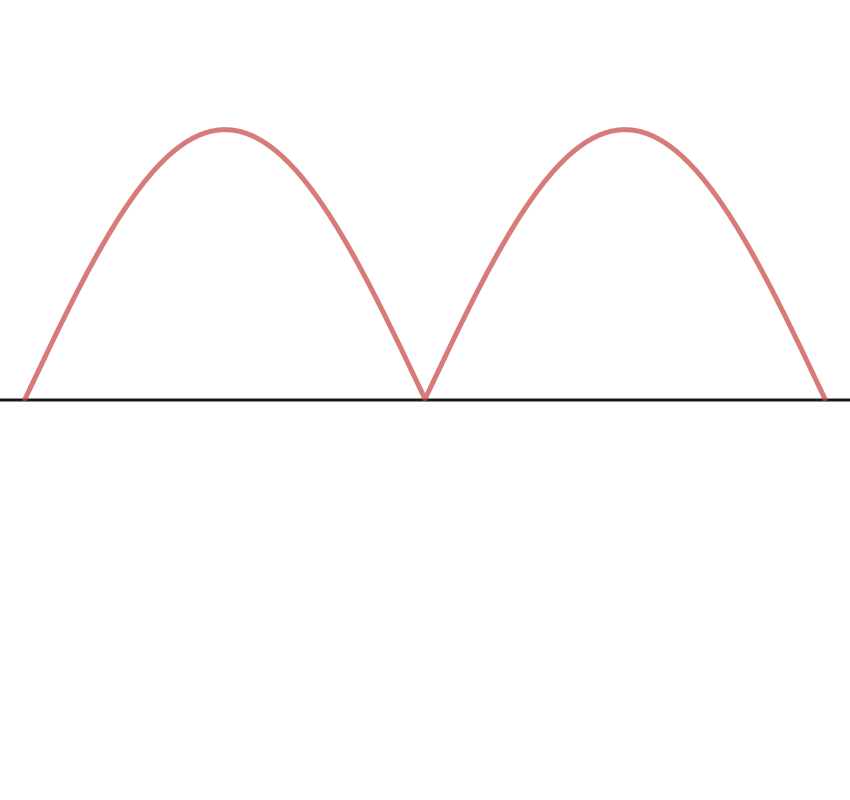}
    \hspace{0.5cm}
    \includegraphics[width=.25\textwidth]{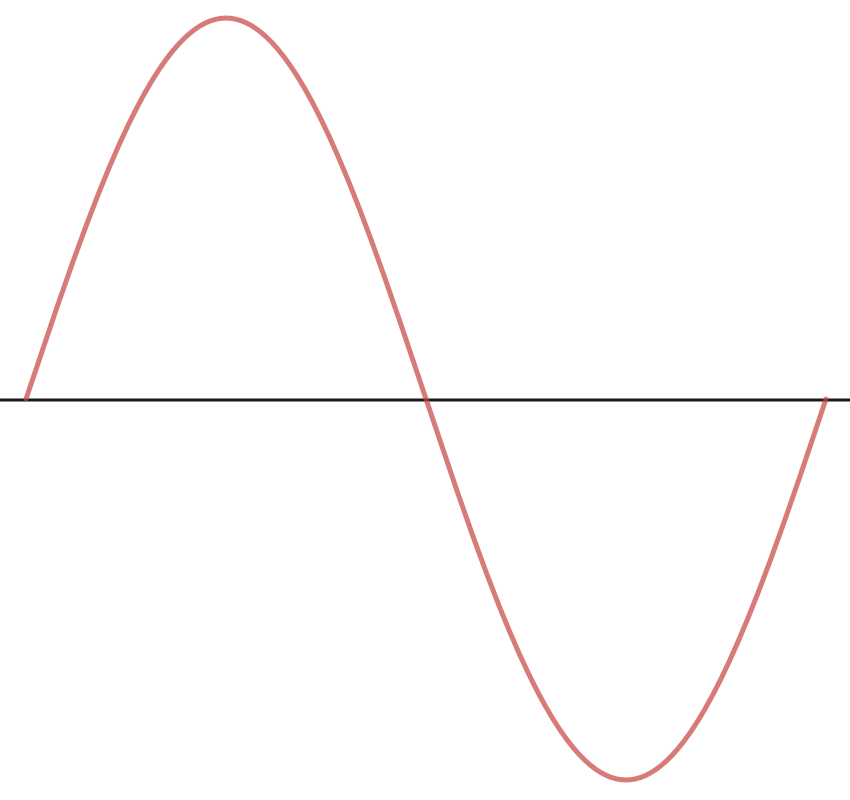}
    \caption{For $n=2$, the limiting eigenfunctions $u_1(x;\infty)$ (left) and $u_2(x;\infty)$ (right)}
    \label{fig1intro}
\end{figure}

\begin{figure}[ht]
  \centering
    \includegraphics[width=0.25\textwidth]{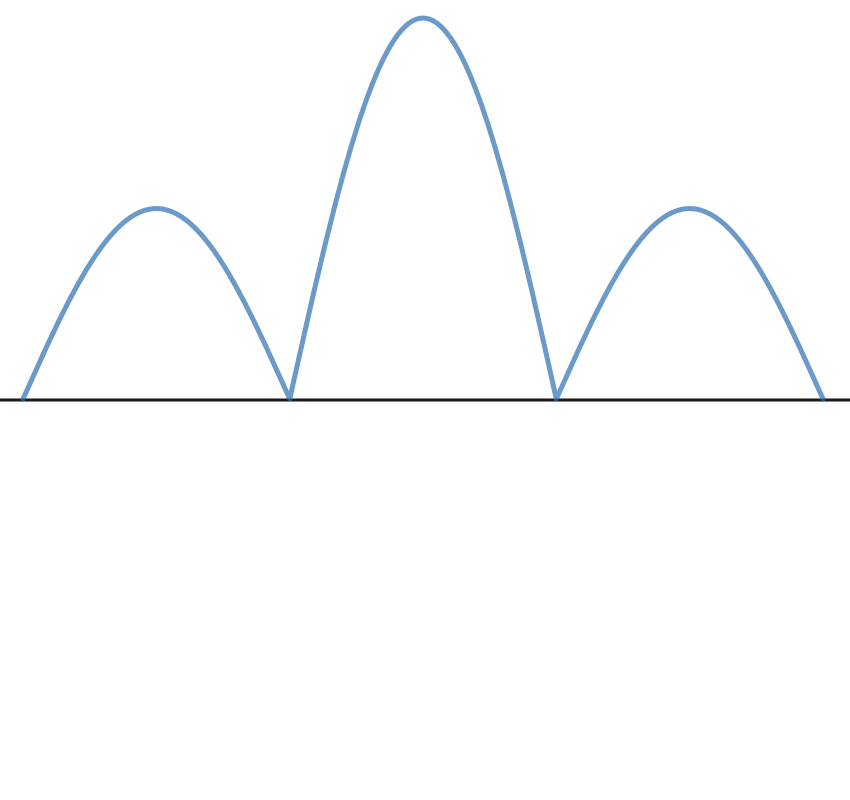}
    \hspace{0.5cm}
    \includegraphics[width=.25\textwidth]{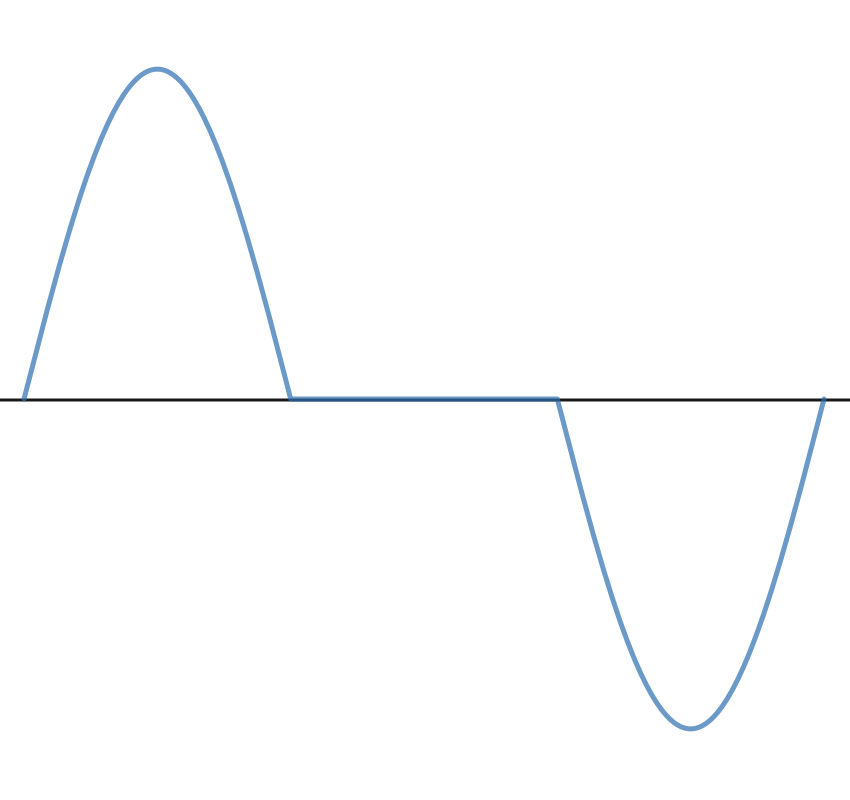}
    \hspace{0.5cm}
    \includegraphics[width=.25\textwidth]{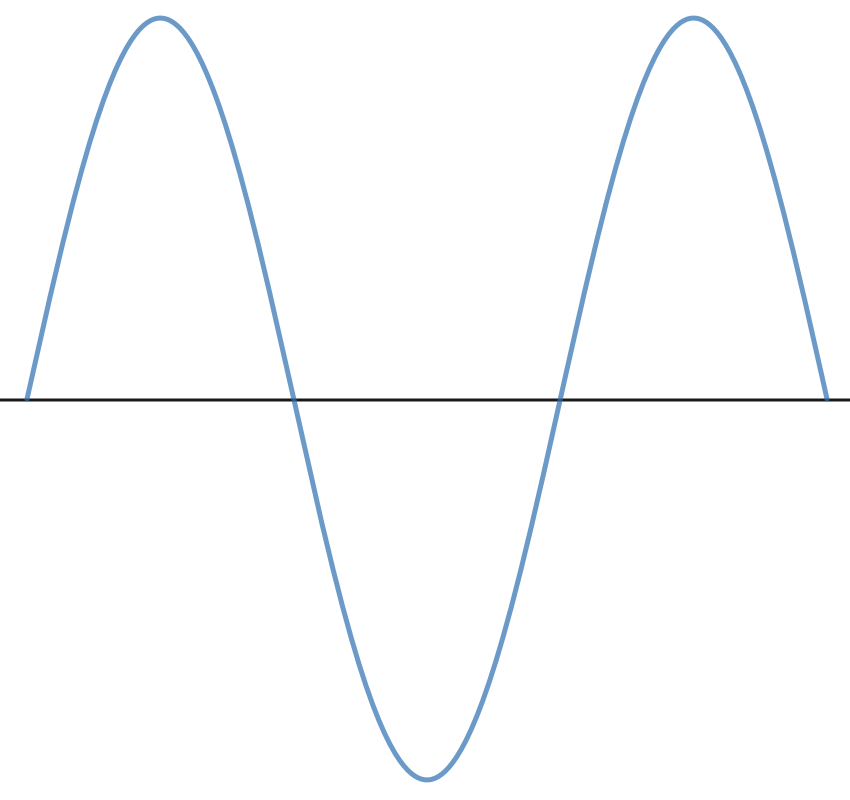}
    \caption{For $n=3$, the limiting eigenfunctions $u_1(x;\infty)$, $u_2(x;\infty)$ and $u_3(x;\infty)$}
    \label{fig2intro}
\end{figure}

For instance, when $n=2$ we have
\begin{align*}
    \big(B_{1,1}, B_{2,1}\big) &= \left( \tfrac{1}{\sqrt2}, -\tfrac{1}{\sqrt2}\right) \\
    \big(B_{1,2}, B_{2,2}\big) &= (1,1)
\end{align*}
corresponding to the left and right sides of Figure \ref{fig1intro}, respectively. Similarly, for $n=3$ we have
\begin{align*}
    \big(B_{1,1}, B_{2,1}, B_{3,1}\big) &= \left( \tfrac12, -1, \tfrac12 \right)\\
    \big(B_{1,2}, B_{2,2}, B_{3,2}\big) &= \left( \tfrac{\sqrt3}2, 0, -\tfrac{\sqrt3}2 \right)\\
    \big(B_{1,3}, B_{2,3}, B_{3,3}\big) &= (1,1,1)
\end{align*}
as shown in Figure \ref{fig2intro}. In general, we see from \eqref{Bkmexplicit} that $B_{\cdot,n} = (1, \ldots, 1)$ for any value of $n$, consistent with the fact that $u_n(x;\sigma) = \sin(n\pi x)$ for all $x \in [0,1]$ and all $\sigma \geq 0$.

\subsection*{Outline of Paper}
We first give the proof of Theorem \ref{thm:main} in Section \ref{sec:proofthm}, and then prove Corollaries \ref{cor:nodes}, \ref{cor:integral} and \ref{cor:infty} in Section \ref{sec:proofcor}. Next, we provide in Section \ref{sec:examples} examples showing that the main theorem and corollaries may fail for non-constant potentials. Lastly, in Section \ref{sec:num} we demonstrate the construction of spectrally accurate numerical methods for the spectral flow, which can be used to verify the results of Theorem \ref{thm:main} for a large range of $n$, and use these methods to demonstrate the behavior of the eigenfunctions for $1\leq n\leq 6$.  We end with a short appendix 
describing our construction by hand of quantum graph versions of the methods in the {\it Chebfun} package developed in \cite{chebfun}.

\subsection*{Acknowledgements} TB was supported in part by NSF Grant DMS-1954304. G.\,Cox acknowledges the support of NSERC grant RGPIN-2017-04259. JLM and G.\,Conte were supported in part by NSF CAREER Grant DMS-1352353 and NSF Applied Math Grant DMS-1909035.  JLM also thanks MSRI for hosting him during the outset of this research project. The authors would like to thank Rowan Killip for asking a question about the eigenfunctions of the spectral flow after a talk by JLM at the Hausdorff Institute in Bonn and hence stimulating part of this research.  The numerics here are an early adaptation of a family of methods being developed by JLM and G.\,Conte in collaboration with Roy Goodman, and the authors thank him for many helpful conversations about these methods.

\section{Proof of Theorem \ref{thm:main}}
\label{sec:proofthm}

To prove the theorem, we first show that the eigenfunctions are given by
\begin{align} \label{eqn:u-m}
 u_m(x;\sigma) = \sin\big(\gamma_m(\sigma)x \big) + \sum_{j=1}^{k-1}A_{j,m}(\sigma) \sin\big(\gamma_m(\sigma)(x- x_j)\big),
\end{align}
for $x\in [x_{k-1},x_k] = \left[\tfrac{k-1}{n},\tfrac{k}{n}\right]$, and with $A_{j,m}(\sigma)$ as in the statement of the theorem. Here we have written $\gamma_{m}(\sigma) = \sqrt{\la_m(\sigma)}$. The function in \eqref{eqn:u-m} clearly satisfies the differential equation
\begin{align*}
u_{m}''(x;\sigma) = - \la_{m}(\sigma)u_{m}(x;\sigma)
\end{align*}
in $(x_{k-1},x_k)$, and so it remains to check the boundary conditions at $x = 0$, $x = 1$, and $x = x_k$. By construction, $u_{m}(0;\sigma) = 0$, and $u_{m}(x;\sigma)$ is continuous at $x = x_k$. Therefore, we are left to choose the coefficients $A_{j,m}(\sigma)$ so that $u_{m}(x;\sigma)$ satisfies the Dirichlet condition at the right end-point $x = 1$, and the derivative jump conditions at $x = x_k$ for each $1 \leq k \leq n-1$. The jump in the derivative at $x = x_k$ is given by
\begin{align*}
u_{m}'(x_k^{+};\sigma) - u_{m}'(x_k^{-};\sigma) = A_{k,m}(\sigma)\gamma_{m}(\sigma),
\end{align*}
and $u_m(x_k;\sigma)$ is equal to
\[
u_{m}(x_k;\sigma) = \sin\left(\tfrac{k\gamma_m(\sigma)}{n} \right) + A_{1,m}(\sigma) \sin\left(\tfrac{(k-1)\gamma_m(\sigma)}{n} \right) + \cdots  + A_{k-1,m}(\sigma)\sin\left(\tfrac{\gamma_m(\sigma)}{n} \right).
\]
Therefore, to ensure that $u_{m}'(x_k^{+};\sigma) - u_{m}'(x_k^{-};\sigma) = \sigma u_{m}(x_k;\sigma)$ and $u_{m}(1;\sigma) = 0$, we require $\gamma_m(\sigma)$ and $A_{j,m}(\sigma)$ to satisfy
\begin{align} \label{eqn:jump1}
 & A_{k,m}(\sigma)\gamma_{m}(\sigma)  = \\
 & \hspace{1cm}  \sigma \left[\sin\left(\tfrac{k\gamma_m(\sigma)}{n} \right) + A_{1,m}(\sigma) \sin\left(\tfrac{(k-1)\gamma_m(\sigma)}{n} \right) + \cdots + A_{k-1,m}(\sigma)\sin\left(\tfrac{\gamma_m(\sigma)}{n} \right)\right] \notag
\end{align}
for $1 \leq k \leq n-1$, and
\begin{align} \label{eqn:Dirichlet}
 \sin\big(\gamma_m(\sigma) \big) + A_{1,m}(\sigma) \sin\left(\tfrac{(n-1)\gamma_m(\sigma)}{n} \right) + \cdots + A_{n-1,m}(\sigma)\sin\left(\tfrac{\gamma_m(\sigma)}{n} \right) = 0.
\end{align}
Our strategy to solve the $n-1$ equations \eqref{eqn:jump1} and equation \eqref{eqn:Dirichlet} is as follows. Using \eqref{eqn:jump1} with $k=1$ to solve for $\sigma$, the remaining $n-2$ equations in \eqref{eqn:jump1} can be written as
\begin{align} \nonumber
&A_{k,m}(\sigma)\sin\left(\tfrac{\gamma_m(\sigma)}{n} \right) = \\ \label{eqn:jump2} & \hspace{1cm} A_{1,m}(\sigma) \left[\sin\left(\tfrac{k\gamma_m(\sigma)}{n} \right) + A_{1,m}(\sigma) \sin\left(\tfrac{(k-1)\gamma_m(\sigma)}{n} \right) + \cdots + A_{k-1,m}(\sigma)\sin\left(\tfrac{\gamma_m(\sigma)}{n} \right)\right] 
\end{align}
for $2\leq k \leq n-1$. Therefore, for fixed $\sigma$, we first will find the $A_{k,m}(\sigma)$ (in terms of $\gamma_{m}(\sigma)$) that solve \eqref{eqn:Dirichlet} and the $n-2$ equations in \eqref{eqn:jump2}. We finally ensure that \eqref{eqn:jump1} holds for $k=1$ by specifying $\gamma_m(\sigma)$ in terms of $\sigma$.
\\
\\
\textbf{Claim:} Given $\gamma_{m}(\sigma)$, the coefficients
\begin{align*}
A_{1,m}(\sigma) & = 2\cos(\theta_m) - 2 \cos\left(\tfrac{\gamma_{m}(\sigma)}{n} \right) , \\
A_{j,m}(\sigma) & = \frac{\sin(j\theta_m)}{\sin(\theta_m)}A_{1,m}(\sigma) \quad\text{ for } 2 \leq j \leq n-1 
\end{align*}
solve the equations \eqref{eqn:Dirichlet} and \eqref{eqn:jump2} for $2 \leq k \leq n-1$. Here $\theta_m = \tfrac{m\pi}{n}$.
\\
\\
\textbf{Proof of Claim:} For ease of notation, we omit the dependence on $\sigma$. We will first prove by induction on $k$ that this choice of $A_{j,m}$ satisfies \eqref{eqn:jump2} for $2 \leq k \leq n-1$. For the base case $k=2$, we have
\begin{align*}
A_{2,m}\sin(\gamma_m/n) = \frac{\sin(2\theta_m)}{\sin(\theta_m)}A_{1,m}\sin(\gamma_m/n) = 2\cos(\theta_m)A_{1,m}\sin(\gamma_m/n) 
\end{align*}
and
\begin{align*}
& A_{1,m} \left[\sin(2\gamma_m/n) + A_{1,m}\sin(\gamma_m/n)\right] \\
& = A_{1,m}\left[\sin(2\gamma_m/n) + 2\cos(\theta_m)\sin(\gamma_m/n)- 2\cos(\gamma_m/n)\sin(\gamma_m/n)\right] \\
& = 2A_{1,m}\cos(\theta_m)\sin(\gamma_m/n) .
\end{align*}
Therefore, \eqref{eqn:jump2} holds for $k=2$. We next assume that \eqref{eqn:jump2} holds for $2 \leq j \leq k$. Then, using the formula for $A_{1,m}$ we have
\begin{align*}
& A_{k,m}\sin(\gamma_m/n) \\
& = A_{1,m}\left[\sin(k\gamma_m/n) + A_{1,m}\sin((k-1)\gamma_m/n) + \cdots + A_{k-1,m}\sin(\gamma_m/n)\right] \\
& = 2\cos(\theta_m)\left[\sin(k\gamma_m/n) + A_{1,m}\sin((k-1)\gamma_m/n) + \cdots + A_{k-1,m}\sin(\gamma_m/n)\right] \\
& - 2 \cos(\gamma_m/n)\left[\sin(k\gamma_m/n) + A_{1,m}\sin((k-1)\gamma_m/n) + \cdots + A_{k-1,m}\sin(\gamma_m/n)\right].
\end{align*}
By the inductive hypothesis (with $j = k$) this implies that
\begin{align*}
& A_{k,m}\sin(\gamma_m/n) =  2\cos(\theta_m)\frac{A_{k,m}}{A_{1,m}}\sin(\gamma_m/n) \\
& \hspace{1cm} -2 \cos(\gamma_m/n)\left[\sin(k\gamma_m/n) + A_{1,m}\sin((k-1)\gamma_m/n) + \cdots + A_{k-1,m}\sin(\gamma_m/n)\right].
\end{align*}
Using
\begin{align*}
2\cos(\gamma_m/n)\sin(j\gamma_m/n) = \sin((j+1)\gamma_m/n) + \sin((j-1)\gamma_m/n),
\end{align*}
we can rewrite the above as
\begin{align*}
& A_{k,m}\sin(\gamma_m/n) =  2\cos(\theta_m)\frac{A_{k,m}}{A_{1,m}}\sin(\gamma_m/n) \\
& \hspace{1.2cm} - \left[\sin((k+1)\gamma_m/n) + A_{1,m}\sin(k\gamma_m/n) + \cdots + A_{k-1,m}\sin(2\gamma_m/n)\right] \\
& \hspace{1.2cm}-  \left[\sin((k-1)\gamma_m/n) + A_{1,m}\sin((k-2)\gamma_m/n) + \cdots + A_{k-2,m}\sin(\gamma_m/n)\right].
\end{align*}
Again by the inductive hypothesis (with $j = k-1$), we can rewrite this equation as
\begin{align}  \nonumber
\sin((k+1)\gamma_m/n) &+ A_{1,m}\sin(k\gamma_m/n) + \cdots + A_{k-1,m}\sin(2\gamma_m/n)+ A_{k,m}\sin(\gamma_m/n) \\ \label{eqn:jump-k+1}
& = 2\frac{A_{k,m}}{A_{1,m}}\cos(\theta_m)\sin(\gamma_m/n) - \frac{A_{k-1,m}}{A_{1,m}}\sin(\gamma_m/n). 
\end{align}
Using the formulae for $A_{j,m}/A_{1,m}$ with $j=k-1$ and $j=k$, the right hand side of \eqref{eqn:jump-k+1} is equal to
\begin{align*}
& \left[2\frac{\sin(k\theta_m)}{\sin(\theta_m)}\cos(\theta_m) - \frac{\sin((k-1)\theta_m)}{\sin(\theta_m)}\cos(\theta_m)\right]\sin(\gamma_m/n) \\
&\qquad  = \frac{\sin((k+1)\theta_m)}{\sin(\theta_m)}\sin(\gamma_m/n) = \frac{A_{k+1,m}}{A_{1,m}}\sin(\gamma_m/n). 
\end{align*}
This completes the proof of the inductive step, and so \eqref{eqn:jump2} holds for all $2 \leq k \leq n-1$. To complete the proof of the claim, we finally need to show that the Dirichlet condition in \eqref{eqn:Dirichlet} is satisfied. To do this, we apply \eqref{eqn:jump-k+1} with $k = n-1$. This gives
\begin{align*}
\sin(\gamma_m) &+ A_{1,m}\sin((n-1)\gamma_m/n) + \cdots + A_{n-1,m}\sin(\gamma_m/n) \\ 
& = 2\frac{A_{n-1,m}}{A_{1,m}}\cos(\theta_m)\sin(\gamma_m/n) - \frac{A_{n-2,m}}{A_{1,m}}\sin(\gamma_m/n) \\
& = \frac{\sin(n\theta_m)}{\sin(\theta_m)}\sin(\gamma_m/n) = \frac{\sin(m\pi)}{\sin(\theta_m)}\sin(\gamma_m/n) = 0
\end{align*}
as required. 
\\
\\
To complete the proof of the theorem, we need to find $\gamma_m(\sigma)$ in terms of $\sigma$ to ensure that \eqref{eqn:jump1} holds for $k=1$, that is
\begin{align*}
A_{1,m}(\sigma)\gamma_m(\sigma) = \sigma \sin\left(\tfrac{\gamma_m(\sigma)}{n}\right).
\end{align*}
Recalling the definition of $A_{1,m}(\sigma)$ and the notation $\gamma_m(\sigma) = \sqrt{\lambda_m(\sigma)}$, this is precisely the implicit equation for $\la_m(\sigma)$ given in the statement of the theorem,
\begin{align} \label{eqn:implicit}
\sigma = 2 \gamma_m(\sigma) \frac{\cos\left(\tfrac{m\pi}{n}\right) - \cos\big(\tfrac{\gamma_m(\sigma)}{n}\big)}{\sin\big(\tfrac{\gamma_m(\sigma)}{n}\big)}.
\end{align}
The right-hand side of \eqref{eqn:implicit} is a strictly increasing function of $\gamma_m(\sigma)$ on the interval $[m\pi,n\pi)$. Moreover, it vanishes when $\gamma_m(\sigma) = m\pi$ and becomes unbounded as $\gamma_m(\sigma)$ approaches $n\pi$, and hence gives a bijection from $[m\pi,n\pi)$ to $[0,\infty)$. Therefore, for each $\sigma\geq0$ there is a unique solution $\gamma_m(\sigma)$. This solution guarantees that \eqref{eqn:jump1} and \eqref{eqn:Dirichlet} both hold, and so gives the desired eigenvalue $\lambda_m(\sigma)$ and eigenfunction $u_m(x;\sigma)$.

\section{Proofs of Corollaries \ref{cor:nodes}, \ref{cor:integral} and \ref{cor:infty}}
\label{sec:proofcor}

\subsection{Proof of Corollary \ref{cor:nodes}}

From our formula for the eigenfunction $u_m(x;\sigma)$ from \eqref{eqn:u-m} we have
\[
u_m(x_k;\sigma) = \sin\left(\tfrac{k\gamma_m(\sigma)}{n} \right) + A_{1,m}(\sigma) \sin\left(\tfrac{(k-1)\gamma_m(\sigma)}{n} \right) + \cdots + A_{k-1,m}(\sigma)\sin\left(\tfrac{\gamma_m(\sigma)}{n} \right).
\]
Using \eqref{eqn:jump2}, this simplifies to
\[
u_m(x_k;\sigma) = \frac{A_{k,m}(\sigma)\sin\big(\tfrac{\gamma_m(\sigma)}{n} \big)}{A_{1,m}(\sigma)}. 
\]
The corollary then follows immediately from the formula for $A_{k,m}(\sigma)/A_{1,m}(\sigma)$ in Theorem \ref{thm:main}.

\subsection{Proof of Corollary \ref{cor:integral}}

Multiplying the eigenfunction equation 
\begin{align*}
    -u_m''(x;\sigma) = \la_m(\sigma)u_m(x;\sigma)
\end{align*}
by $\sin(n\pi x)$ and integrating by parts twice, we have
\begin{align*}
-\la_m(\sigma)\int_{I_k}u_{m}(x;\sigma) \sin(n\pi x) \, dx  &= \int_{I_k}u_{m}''(x;\sigma) \sin(n\pi x) \, dx \\
& = \left[u_{m}'(x;\sigma)\sin(n\pi x) - n\pi u_{m}(x;\sigma)\cos(n\pi x) \right]_{x_{k-1}}^{x_k} \\
& \hspace{0.5cm} -n^2\pi^2 \int_{I_k}u_{m}(x;\sigma) \sin(n\pi x) \, dx.
\end{align*}
Since $\sin(n\pi x)$ vanishes at the nodes $x_{k-1}$ and $x_{k}$, while $\cos(n\pi x_j) = (-1)^{j}$, this simplifies to
\begin{align*}
-\la_m(\sigma)F_{k,m}(\sigma) = -n\pi(-1)^ku_{m}(x_k;\sigma) + n\pi (-1)^{k-1}u_{m}(x_{k-1};\sigma) - n^2\pi^2F_{k,m}(\sigma). 
\end{align*}
Rearranging this gives the equality in the corollary. Moreover, if $F_{1,m}(\sigma) \equiv 1$, then
\begin{align*}
F_{k,m}(\sigma) = (-1)^{k+1}\frac{u_m(x_k;\sigma) + u_{m}(x_{k-1};\sigma)}{u_{m}(x_1;\sigma)}.
\end{align*}
By Corollary \ref{cor:nodes}, this simplifies to
\begin{align*}
F_{k,m}(\sigma) = (-1)^{k+1}\frac{u_m(x_k;0) + u_{m}(x_{k-1};0)}{u_{m}(x_1;0)},
\end{align*}
which is independent of $\sigma$.

\subsection{Proof of Corollary \ref{cor:infty}}

At $\sigma = 0$, we have the eigenfunction $u_m(x;0) = \sin(m\pi x)$, and this leads to
\begin{align*}
F_{k,m}(0) & = \int_{I_k} \sin(m \pi x)\sin(n \pi x) \, dx \\
& = (-1)^{k+1} \frac{2n}{\pi(n-m)(n+m)} \sin{\frac{(2k-1)m\pi}{2n}} \cos{\frac{m\pi}{2n}}.
\end{align*}
Using Corollary \ref{cor:integral}, we choose a normalization of $u_m(x;\sigma)$ for $\sigma >0$ so that $F_{k,m}(\sigma)$ is constant in $\sigma$. As $\sigma$ tends to $\infty$, $u_m(x;\sigma)$ tends to $B_{k,m}\sin(n\pi x)$ on $I_k$ for some $B_{k,m}$. This constant $B_{k,m}$ must therefore satisfy
\begin{align*}
F_{k,m}(0) = \lim_{\sigma\to\infty}F_{k,m}(\sigma) = B_{k,m}\int_{I_k}\sin^2(n\pi x) \, dx.
\end{align*}
Solving for $B_{k,m}$ and rescaling by a factor that does not depend on $k$ proves the corollary.

\section{Failure under modifications of the assumptions} \label{sec:examples}

In this section we demonstrate that the results in Theorem \ref{thm:main} and Corollaries \ref{cor:nodes}, \ref{cor:integral} and \ref{cor:infty} are particular to the spectral flow for the Laplacian (i.e. the case $V \equiv 0$) on an interval. For a general potential $V$, we can write
\begin{align*}
F_{k,m}^{V}(\sigma) = \int_{I_k}u_m(x;\sigma)u_n(x) \, dx, 
\end{align*}
where $u_n$ is the $n$-th eigenfunction of \eqref{Vdiffeq} and $u_m(x;\sigma)$ the $m$-th eigenfunction of the associated spectral flow in \eqref{Vdiffsum}. In this case, we integrate over the interval $I_k = [x_{k-1},x_k]$, where $x_k$ are the zeros of $u_n$. Following the integration by parts calculation in the proof of Corollary \ref{cor:integral}, we obtain
\begin{align*}
F_{k,m}^{V}(\sigma) = \frac{-u_m(x_k;\sigma)u_n'(x_k) + u_{m}(x_{k-1};\sigma)u_n'(x_{k-1})}{\la_n - \la_m(\sigma)}.
\end{align*}
However, unlike the case $V=0$, for non-zero potentials the quantities $F_{k,m}^{V}(\sigma)$ will in general not be independent of $\sigma$ for any normalization of $u_m(x;\sigma)$. This is because $u_m(x_k;\sigma)/u_m(x_1;\sigma)$ is not in general independent of $\sigma$ for $V\neq0$.

We will prove this for a generic small perturbation of the potential from $V=0$, and also give a numerical demonstration below.

Let $V = \eps W$, where $W$ is a continuous function on $[0,1]$, and $\eps>0$ is a small constant. Denoting by $u_{n,\eps}(x)$ the eigenfunctions of 
\begin{align*}
-u''(x) + \eps Wu(x) = \la u(x), \qquad u(0) = u(1) = 0
\end{align*}
we set $n=3$, and let $x_1(\eps)$, $x_2(\eps)$ be the interior zeros of $u_{3,\eps}(x)$. We then consider the spectral flow
\begin{align*}
-u''(x) + \eps Wu(x) + \sigma \big[\delta\big(x-x_1(\eps)\big) + \delta\big(x-x_2(\eps)\big)\big]u(x) = 0, \qquad u(0) = u(1) = 0.
\end{align*}
Letting $\la_m(\eps,\sigma)$ and $u_m(x;\eps,\sigma)$ denote the eigenvalues and eigenfunctions of this spectral flow, with $m=1,2$, we set
\begin{align} \nonumber
F_{k,m}(\eps,\sigma) & = \int_{I_k(\eps)}u_m(x;\eps,\sigma)u_{3,\eps}(x) \, dx \\  \label{eqn:FkV}
& = \frac{-u_m(x_k(\eps);\eps,\sigma)u_{3,\eps}'(x_k(\eps)) + u_{m}(x_{k-1}(\eps);\eps,\sigma)u_{3,\eps}'(x_{k-1}(\eps))}{\la_3(\eps,0) - \la_m(\eps,\sigma)}.
\end{align}
Here $I_k(\eps) = [x_{k-1}(\eps),x_k(\eps)]$ for $1\leq k\leq3$, where we have set $x_0(\eps)=0$, $x_3(\eps) = 1$.  When $\eps=0$, all these quantities reduce to those in Theorem \ref{thm:main} and Corollary \ref{cor:integral}.

 \begin{theorem} \label{thm:counter}
 Suppose that
\begin{align*}
    c_W := \int_{0}^{1}W(t)\sin(m\pi(1-t))\, dt \neq 0.
\end{align*}
 Then, there exists a constant $\eps_0 = \eps_0\left(c_W,\norm{W}_{L^{\infty}}\right) >0$ such that for all $0<\eps<\eps_0$, no normalization of the eigenfunctions $u_m(x;\eps,\sigma)$ can ensure that the quantities $F_{k,m}(\eps,\sigma)$  are simultaneously independent of $\sigma$ for $1\leq k \leq 3$.
 \end{theorem}
  \begin{proof1}{Theorem \ref{thm:counter}}
Suppose that for $\eps>0$ fixed we choose a normalization of the eigenfunctions $u_m(x;\eps,\sigma)$ for $\sigma \geq0$ so that $F_{1,m}(\eps,\sigma)$ is independent of $\sigma$. Then, from \eqref{eqn:FkV}, $F_{2,m}(\eps,\sigma)$ will be independent of $\sigma$ precisely when the quantity $u_m(x_2;\eps,\sigma)/u_m(x_1;\eps,\sigma)$ is independent of $\sigma$. We show that this does not hold for $\eps$ and $\sigma$ sufficiently small. To see this, we first note that the equivalent ansatz to \eqref{eqn:u-m} in the proof of Theorem \ref{thm:main} is to write
\begin{align*}
u_m(x;\eps,\sigma) = w_{0,m}(x) + \sum_{j=1}^{k-1}A_{j,m}(\eps,\sigma)w_{j,m}(x-x_j(\eps)) \end{align*} 
for $x\in I_k(\eps) = [x_{k-1}(\eps),x_k(\eps)]$.  Here the functions $w_{k,m}(x;\eps,\sigma)$ satisfy
\begin{align} \label{eqn:wk}
-w_{k,m}''(x;\eps,\sigma) +\eps W(x+x_k(\eps))w_{k,m}(x;\eps,\sigma)  = \la_m(\eps,\sigma)&w_{k,m}(x;\eps,\sigma), \\\nonumber w_{k,m}(0;\eps,\sigma) & = 0, \quad w_{k,m}'(0;\eps,\sigma) = \gamma_m(\eps,\sigma)
\end{align} 
with $\gamma_m(\eps,\sigma) = \sqrt{\la_m(\eps,\sigma)}$. In particular, we have
\begin{align*}
\left|x_1(\eps) -\tfrac{1}{3}\right| + \left|x_2(\eps) -\tfrac{2}{3}\right| + \big|w_{k,m}(x;\eps,\sigma) - \sin\left(\gamma_{m}(\sigma)x\right)\big| = O(\eps).
\end{align*}
We use the $O$ notation to denote a constant depending only on the $L^{\infty}$-norm of $W$ (but independent of $\eps$ and $\sigma$). Suppose for contradiction that $u_m(x_1(\eps);\eps,\sigma)/u_m(x_2(\eps);\eps,\sigma)$ is independent of $\sigma$. Then, $A_{1,m}(\eps,\sigma)$ must satisfy
\begin{align*}
A_{1,m}(\eps,\sigma)\frac{w_{1,m}(x_2(\eps)-x_1(\eps);\eps,\sigma)}{w_{0,m}(x_1(\eps);\eps,\sigma)} + \frac{w_{0,m}(x_2(\eps);\eps,\sigma)}{w_{0,m}(x_1(\eps);\eps,\sigma)} = \frac{w_{0,m}(x_2(\eps);\eps,0)}{w_{0,m}(x_1(\eps);\eps,0)}.
\end{align*} 
In particular, $A_{1,m}(\eps,\sigma) = A_{1,m}(\sigma) + O(\eps \sigma)$, with $A_{1,m}(\sigma) = A_{1,m}(0,\sigma) = O(\sigma)$ as in the statement of Theorem \ref{thm:main}. To ensure that the jump conditions at $x=x_1(\eps)$ and $x=x_2(\eps)$ hold, analogously to \eqref{eqn:jump1},  $A_{2,m}(\eps,\sigma)$ must then satisfy
\begin{align*}
A_{2,m}(\eps,\sigma)w_{0,m}(x_1(\eps);\eps,\sigma) = A_{1,m}(\eps,\sigma)\left[w_{0,m}(x_2(\eps);\eps,\sigma) + A_{1,m}(\eps,\sigma)w_{1,m}(x_2(\eps);\eps,\sigma) \right].
\end{align*}
This again ensures that $A_{2,m}(\eps,\sigma) = A_{2,m}(\sigma) + O(\eps\sigma)$. Finally, to ensure that the Dirichlet condition $u_{m}(1;\eps,\sigma) = 0$ holds, we require that
\begin{align} \label{eqn:eps1}
w_{0,m}(1;\eps,\sigma) + A_{1,m}(\eps,\sigma)w_{1,m}(1-x_1(\eps)) + A_{2,m}(\eps,\sigma)w_{2,m}(1-x_2(\eps)) = 0.
\end{align}
From Theorem \ref{thm:main} we know that \eqref{eqn:eps1} holds when $\eps = 0$, and so we can rewrite it as
\begin{align} \label{eqn:eps2}
w_{0,m}(1;\eps,\sigma) - w_{0,m}(1;0,\sigma) + O(\eps \sigma) = 0.
\end{align}
Also, from \eqref{eqn:wk} the function $w_{0,m}(x;\eps,\sigma)$ satisfies
\begin{align*}
w_{0,m}(1;\eps,\sigma) & =  \sin\left(\gamma_{m}(\sigma)(1)\right) +\eps\int_{0}^{1}W(t)\sin(\gamma_m(0)(1-t))\, dt + O(\eps \sigma+\eps^2) \\
& =  w_{0,m}(1;0,\sigma) +\eps\int_{0}^{1}W(t)\sin(m\pi(1-t))\, dt + O(\eps \sigma+\eps^2).
\end{align*}
Therefore, \eqref{eqn:eps2} becomes
\begin{align*}
\eps\int_{0}^{1}W(t)\sin(m\pi(1-t))\, dt + O(\eps\sigma+\eps^2) = 0.
\end{align*}
Since it was assumed that this integral does not vanish, we obtain at a contradiction for all $\eps$ and $\sigma$ sufficiently small. Therefore, $u_m(x_2;\eps,\sigma)/u_m(x_1;\eps,\sigma)$ is not independent of $\sigma$, and this ensures that $F_{k,m}(\eps,\sigma)$ cannot be independent of $\sigma$ for $k=1$ and $k=2$ simultaneously.
\end{proof1}

\begin{figure}[ht]
  \centering
    \includegraphics[width=0.45\textwidth]{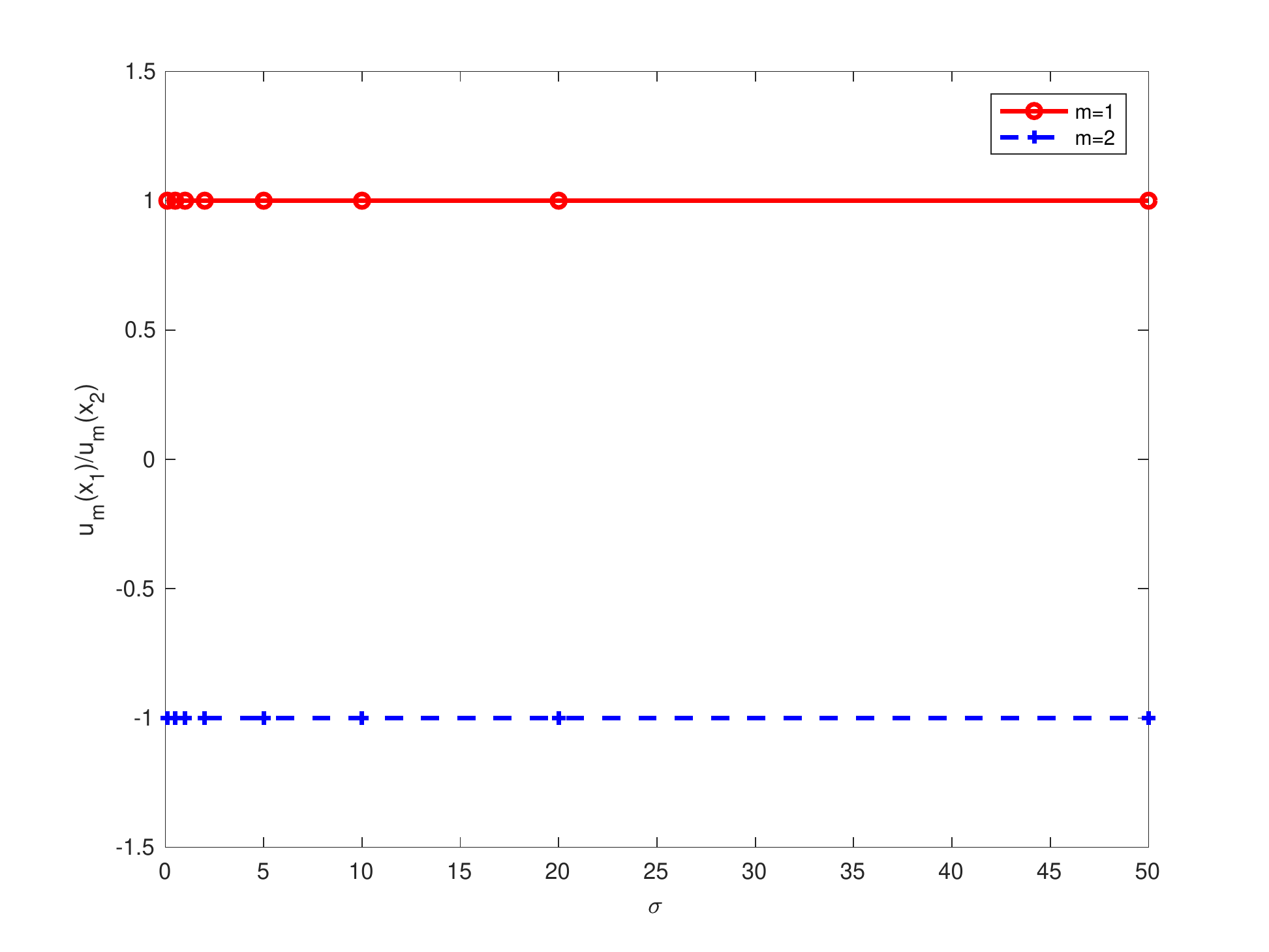}
    \includegraphics[width=.45\textwidth]{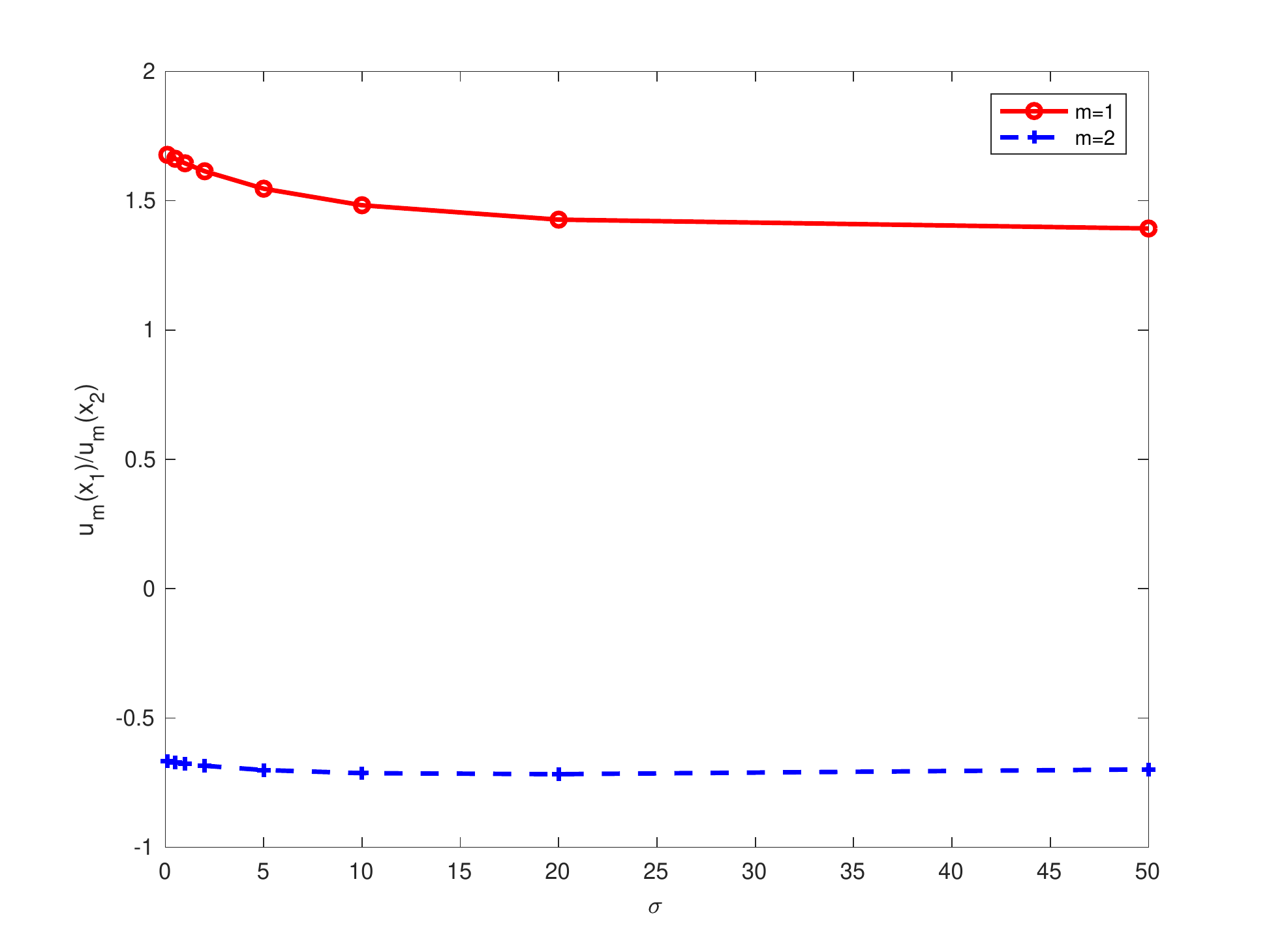}
    \caption{The ratio of $u_m(x_1;\sigma)/u_m(x_2;\sigma)$ for the spectral flow acting on eigenfunctions of the operators $H_0 = -\partial_x^2$ (left) and $H_1 = -\partial_x^2+ 20\rchi_{(0,1/2)}$ (right)}
    \label{fail}
\end{figure}

We now illustrate the dependence of the quantities $F_{k,m}^V(\sigma)$ on $\sigma$ numerically. Using computational methods described in \cite{roy}, we numerically compute the spectrum and the corresponding eigenfunctions the operators $H_0 = -\partial_x^2$ and $H_1 = -\partial_x^2 + 20 \rchi_{(0,1/2)}$ on $[0,1]$ where $\rchi_{(0,1/2)}$ is the characteristic function that takes values $1$ for $0 < x < 1/2$ and $0$ otherwise.  By computing the third eigenfunction of these operators with a very densely defined grid ($768$ points on the interval) and using linear interpolation to compute the approximate zeros and set the spectral flow boundary conditions accordingly, we  construct the corresponding spectral flow based upon the computed nodal set.  The results are plotted in Figure \ref{fail}, where it is shown that for the spectral flow generated by the nodal set of the third eigenfunction, the ratios of $u_m(x_1;\sigma)/u_m(x_2;\sigma)$ for $m=1,2$ are independent of $\sigma$ in our computation for $H_0$, but depend upon $\sigma$ in a nonlinear (and in particular nonconstant) fashion for $H_1$.
\\

\section{Numerical Methods and Results}
\label{sec:num}

We use a numerical method in MATLAB to approximate the eigenfunctions of \eqref{diffeq} and \eqref{deltabc}. In the following section we describe the method in detail and  display results for different values of $n$.

\subsection{Chebyshev Discretization Method}

The method used to approximate the eigenfunctions is based on a Chebyshev discretization using rectangular differentiation matrices \cite{xu}. Traditionally, spectral collocation methods have involved the deletion and replacement of rows of a square matrix in order to impose boundary conditions; the rectangular collocation method employed here instead uses a resampling of the interpolating polynomials to produce rectangular matrices without needing to delete any rows \cite{hale}. 

In our method, the differentiation matrix is applied on a quantum graph. A quantum graph is a metric graph, i.e. a set of vertices and edges, where each edge connects two vertices and is assigned a positive length. An operator can then be defined on the edges, and boundary conditions imposed at the vertices (see \cite{graph} for an introduction to quantum graphs). Here our graph has $n$ edges connecting $n+1$ vertices (nodes). We first build a quantum graph with the desired number of nodes, subinterval lengths, and boundary conditions, then apply the operator and solve the eigenvalue problem on the graph using a Chebychev generalization of a quantum graph package developed by R. Goodman; more on this package can be found in \cite{roy}.

Finally, we write a function to obtain a matrix containing the amplitudes of each eigenfunction on each subinterval, normalize this matrix, and find the maximum difference between the values of this matrix and the matrix of the eigenvectors from Corollary \ref{cor:infty}. We also find the difference between the first and $n$-th eigenvalues for each $n$. If small enough, these values will imply the convergence of the eigenvalues and an agreement with Theorem \ref{thm:main}.

\subsection{Results}

For small values of $n$ we display the progression of eigenvalues as $\sigma$ increases, as well as plots of the eigenfunctions when $\sigma = 10^7$. The amplitudes of the eigenfunctions at $\sigma =\infty$ from Corollary \ref{cor:infty}
are given by the matrix $M_{\rm thm}$. The entry $M_{\rm thm}[i,j]$ represents the coefficient of $\sin({n \pi x})$ on the $j$-th sub-interval of the $i$-th eigenfunction, with $1 \leq i \leq n-1$ and $1 \leq j \leq n$. The matrix $M_{\rm norm}$ consists of the normalized vectors produced using the Chebyshev method. We give the maximum difference between the numerically observed and actual matrices $M_{\rm norm}$ and $M_{\rm thm}$ for each $n$, which we will refer to as diff$_{\rm vec}$, and the difference between the first and $n$-th eigenvalues, which we will call diff$_{\rm val}$.  \\

\paragraph{($n = 2$)}
The following table  shows the first two eigenvalues, $\lambda_1$ and $\lambda_2$, with increasing $\sigma$:

\begin{center}
\begin{tabular}{|p{1cm}|p{1.5cm}|p{1.5cm}|p{1.5cm}|p{1.5cm}|p{1.5cm}| }
 \hline
  & $\sigma = 0$ & $\sigma = 10$ & $\sigma = 10^3$ & $\sigma = 10^5$ & $\sigma = 10^7$ \\
 \hline
 $\lambda_1$ & 9.8696 & 22.6699 & 39.1645 & 39.4753 & 39.4784\\
 $\lambda_2$ & 39.4784 & 39.4784 & 39.4784 & 39.4784 & 39.4784\\
 \hline
\end{tabular}
\end{center}

\begin{figure}[ht]
  \centering
    \includegraphics[width=.3\textwidth]{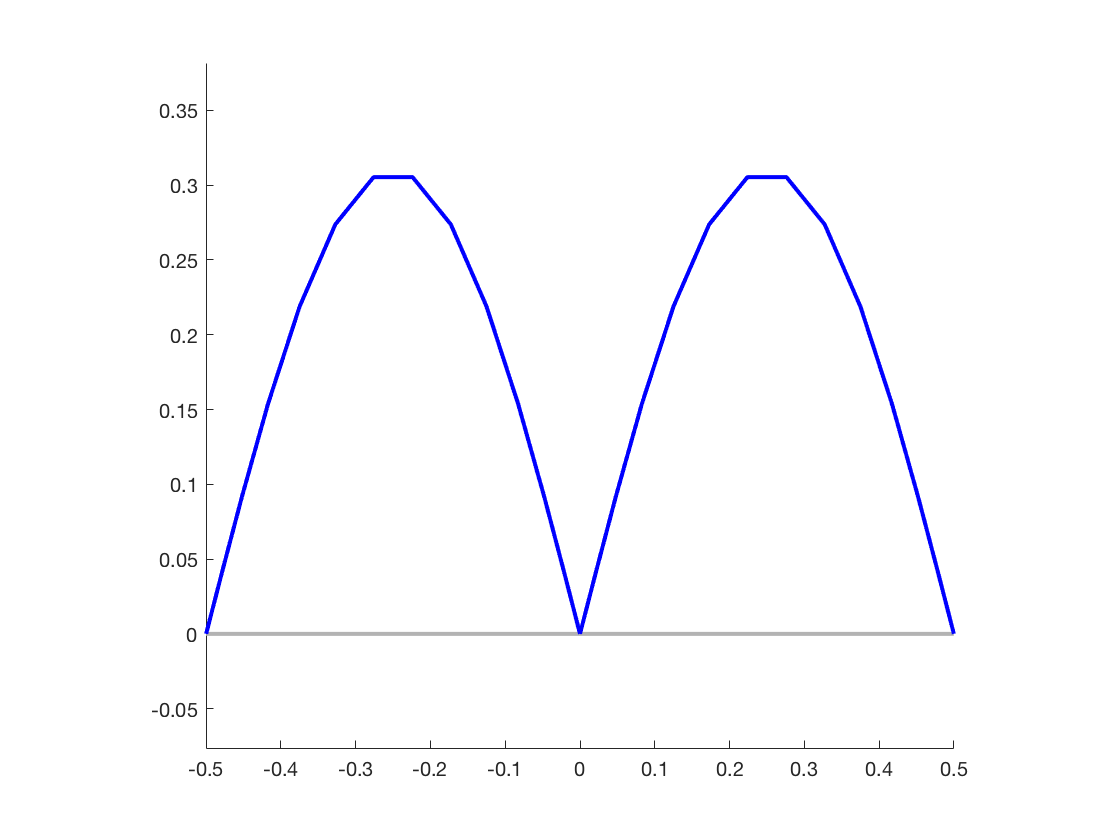}
    \includegraphics[width=.3\textwidth]{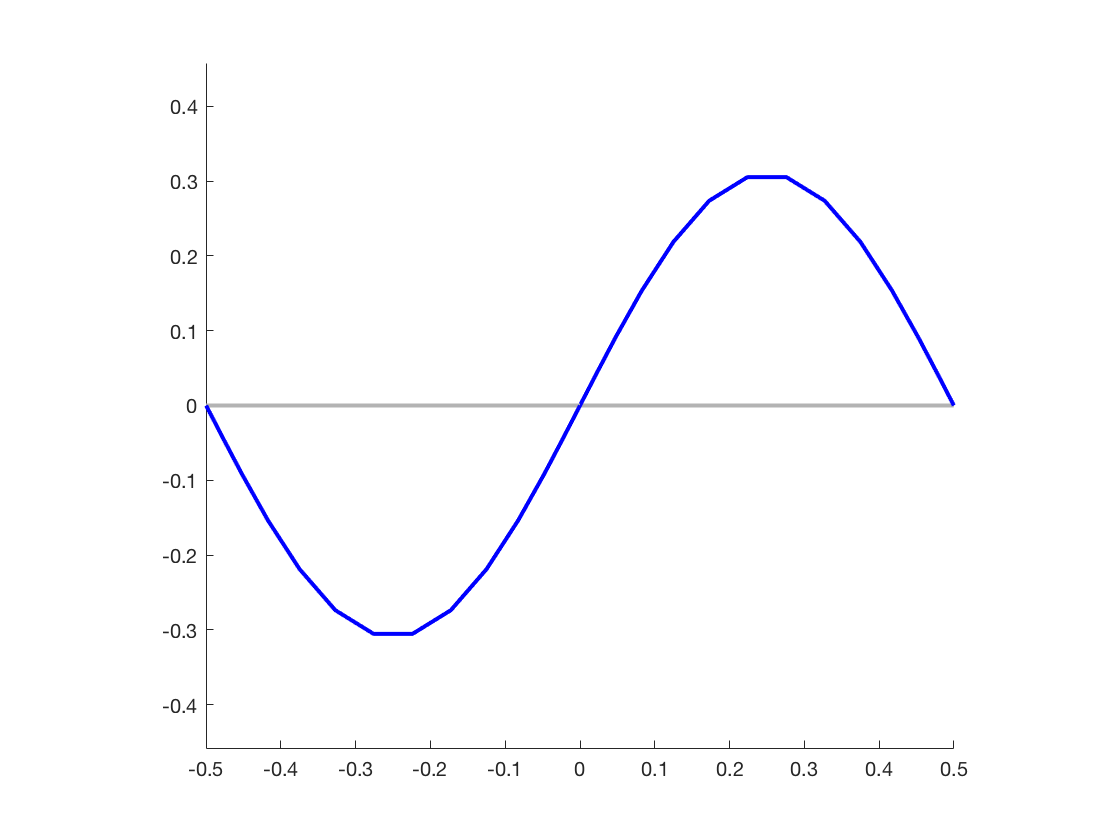}
    \caption{(Left) $u_{1,2}(x)$, $\lambda_1=39.4784$.  (Right) $u_{2,2}(x)$, $\lambda_2=39.4784$.}
    \label{fig1}
\end{figure}

Observe that $\lambda_1(\sigma)$ is strictly increasing while 
$\lambda_2(\sigma)$ is constant, as expected. The vector corresponding to the first eigenfunction is given by $M_{\rm thm} = (0.7071 ,-0.7071)$, 
while the normalized vector obtained using the Chebyshev method with $\sigma = 10^7$ is $M_{\rm norm}= (0.7071 ,-0.7071)$.
Based on these two matrices, we obtain the following for the values of diff$_{\rm vec}$ and diff$_{\rm val}$:

\begin{center}
\begin{tabular}{ |p{3cm}|p{3cm}| }
 \hline
{diff}$_{\rm val}$ & {diff}$_{\rm vec}$ \\
 \hline
$ 8.0000 * 10^{-7}$ & $1.6106 * 10^{-8}$ \\
 \hline
\end{tabular}
\end{center}
The eigenfunctions $u_{1,2}(x;\sigma)$ and $u_{2,2}(x;\sigma)$ at $\sigma = 10^7$ are displayed in Figure \ref{fig1}. \\

\paragraph{($n = 3$)}

For the case of three subintervals, the first three eigenvalues
are given in the following table: 

\begin{center}
\begin{tabular}{|p{1cm}|p{1.5cm}|p{1.5cm}|p{1.5cm}|p{1.5cm}|p{1.5cm}| }
 \hline
  & $\sigma = 0$ & $\sigma = 10$ & $\sigma = 10^3$ & $\sigma = 10^5$ & $\sigma = 10^7$ \\
 \hline
 $\lambda_1$ & 9.8696 & 32.6297 & 87.2491 & 88.8105 & 88.8263\\
 $\lambda_2$ & 39.4784 & 59.8161 & 88.2959 & 88.8211 & 88.8264\\
 $\lambda_3$ & 88.8264 & 88.8264 & 88.8264 & 88.8264 & 88.8264 \\
 \hline
\end{tabular}
\end{center}

\begin{figure}[ht]
  \centering
    \includegraphics[width=0.3\textwidth]{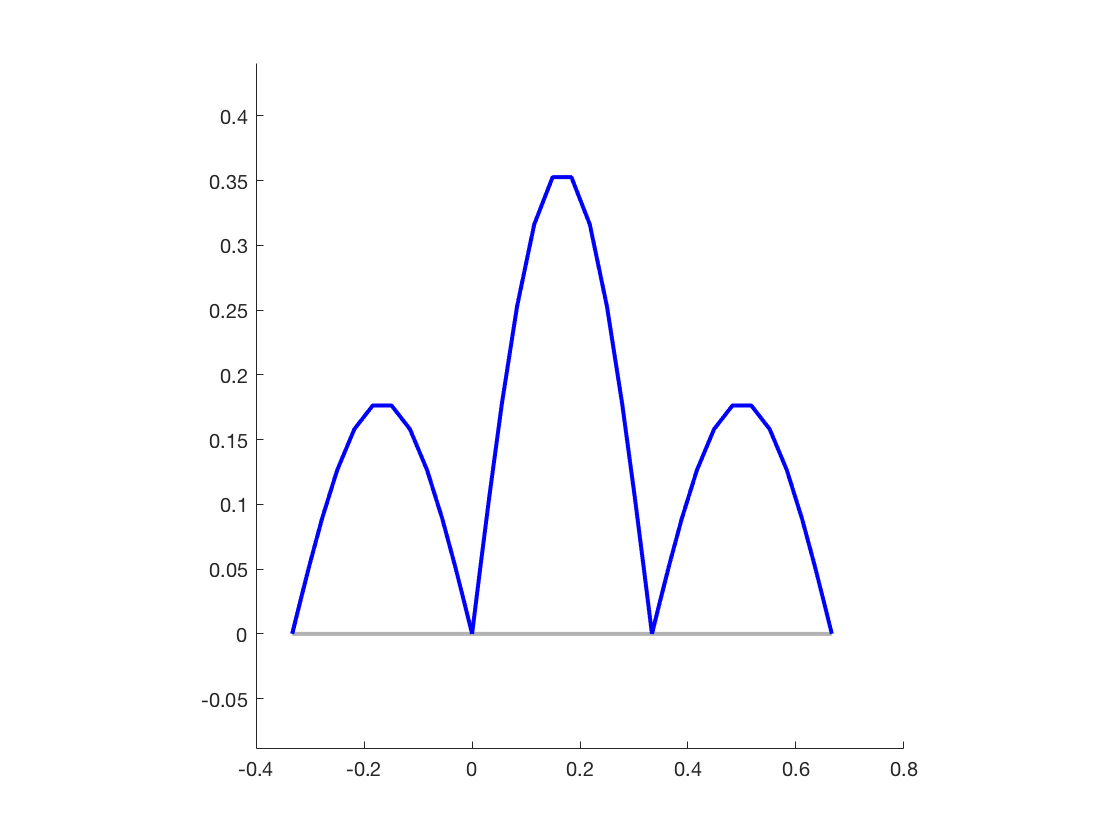}
    \includegraphics[width=.3\textwidth]{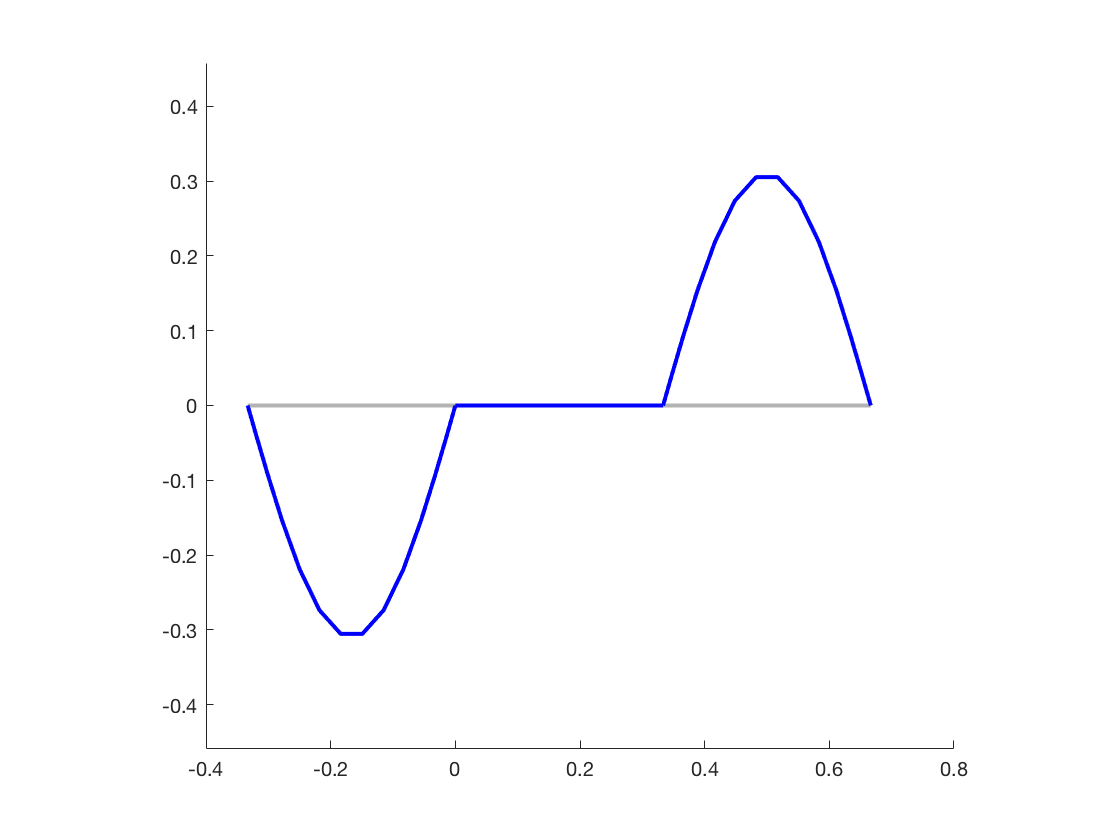}
    \includegraphics[width=.3\textwidth]{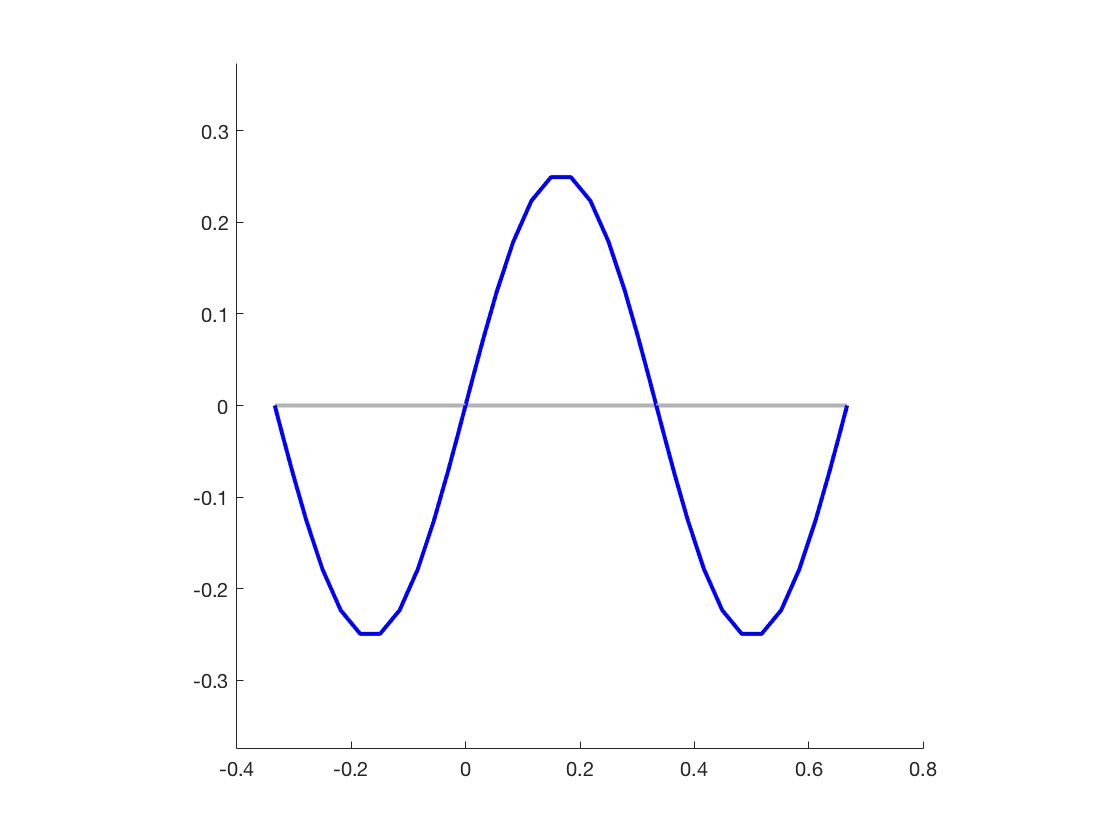}
    \caption{(Left) $u_{1,3}(x)$, $\lambda_1=88.8263$.  (Center) $u_{2,3}(x)$, $\lambda_2=88.8264$. (Right) $u_{3,3}(x)$, $\lambda_3=88.8264$.}
    \label{fig2}
\end{figure}

As expected, we see that $\lambda_1(\sigma)$ and $\lambda_2(\sigma)$ are strictly increasing and $\lambda_3(\sigma)$ is constant.
The matrices $M_{\rm thm}$ and $M_{\rm norm}$ for $n=3$ with $\sigma = 10^7$ are then
\begin{equation*}
    M_{\rm thm} = 
\begin{pmatrix}
0.5000  & -1.0000 & 0.5000 \\
0.8660 & -0.0000 & -0.8660 \\
\end{pmatrix}
\end{equation*}
and
\begin{equation*}
    M_{\rm norm} = 
\begin{pmatrix}
0.5000  & -1.0000 & 0.5000 \\
0.8660 & 0.0000 & -0.8660 \\
\end{pmatrix}.
\end{equation*}
From these matrices, we obtain the following:
\begin{center}
\begin{tabular}{ |p{3cm}|p{3cm}| }
 \hline
 diff$_{\rm val}$ & {diff}$_{\rm vec}$ \\
 \hline
 $1.8000 * 10^{-6}$ & $8.2734 * 10^{-7} $\\
 \hline
\end{tabular}
\end{center}
The plots of the first three eigenfunctions $u_{1,3}(x)$, $u_{2,3}(x)$, and $u_{3,3}(x)$ at $\sigma = 10^7$ are displayed in Figure \ref{fig2}. \\

\paragraph{($n = 4$)}

For $n=4$, we obtain the following eigenvalues:

\begin{center}
\begin{tabular}{|p{1cm}|p{1.5cm}|p{1.5cm}|p{1.5cm}|p{1.5cm}|p{1.5cm}| }
 \hline
  & $\sigma = 0$ & $\sigma = 10$ & $\sigma = 10^3$ & $\sigma = 10^5$ & $\sigma = 10^7$ \\
 \hline
 $\lambda_1$ & 9.8696 & 42.4846 & 153.6882 & 157.8705 & 157.9132\\
 $\lambda_2$ & 39.4784 & 70.9891 & 155.4176  & 157.8884 & 157.9134\\
 $\lambda_3$ & 88.8264 & 1116.7243 & 157.1763 & 157.9063 & 157.9136 \\
 $\lambda_4$ & 157.9137 & 157.9137 & 157.9137 & 157.9137 &  157.9137 \\
 \hline
\end{tabular}
\end{center}

\begin{figure}[ht]
  \centering
    \includegraphics[width=.3\textwidth]{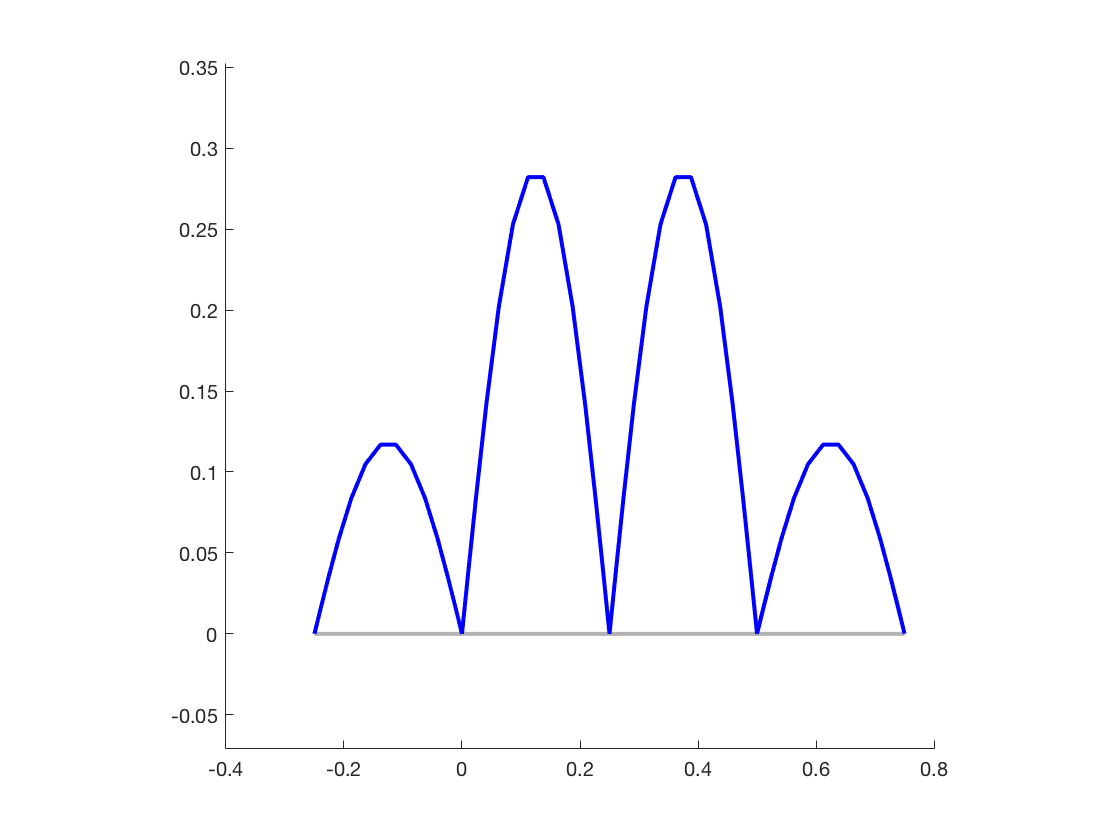}
    \includegraphics[width=.3\textwidth]{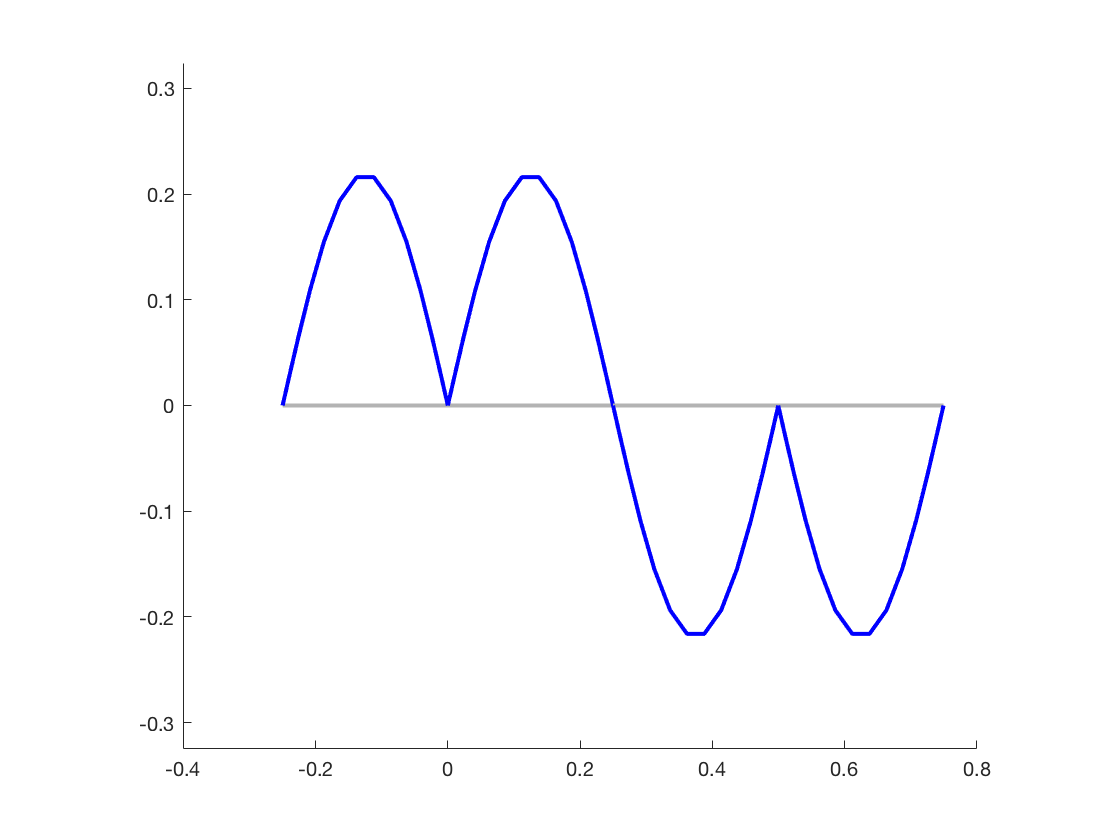} \\
    \includegraphics[width=.3\textwidth]{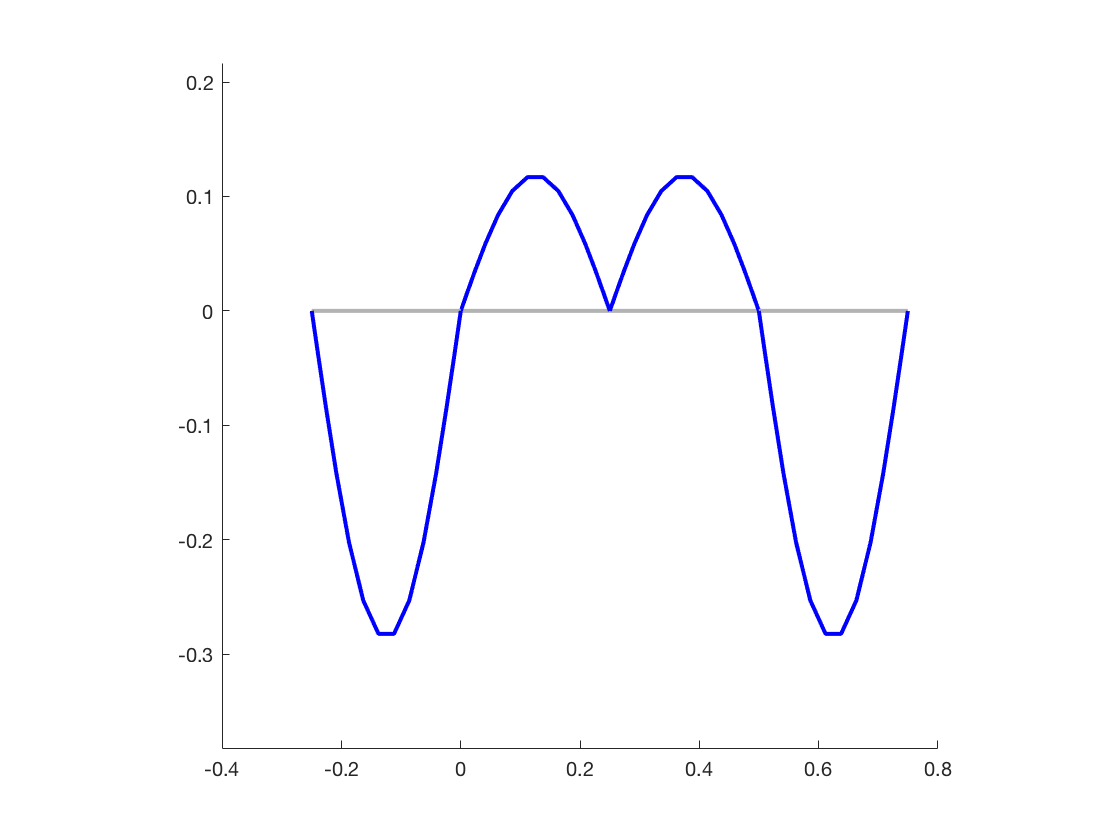}
    \includegraphics[width=.3\textwidth]{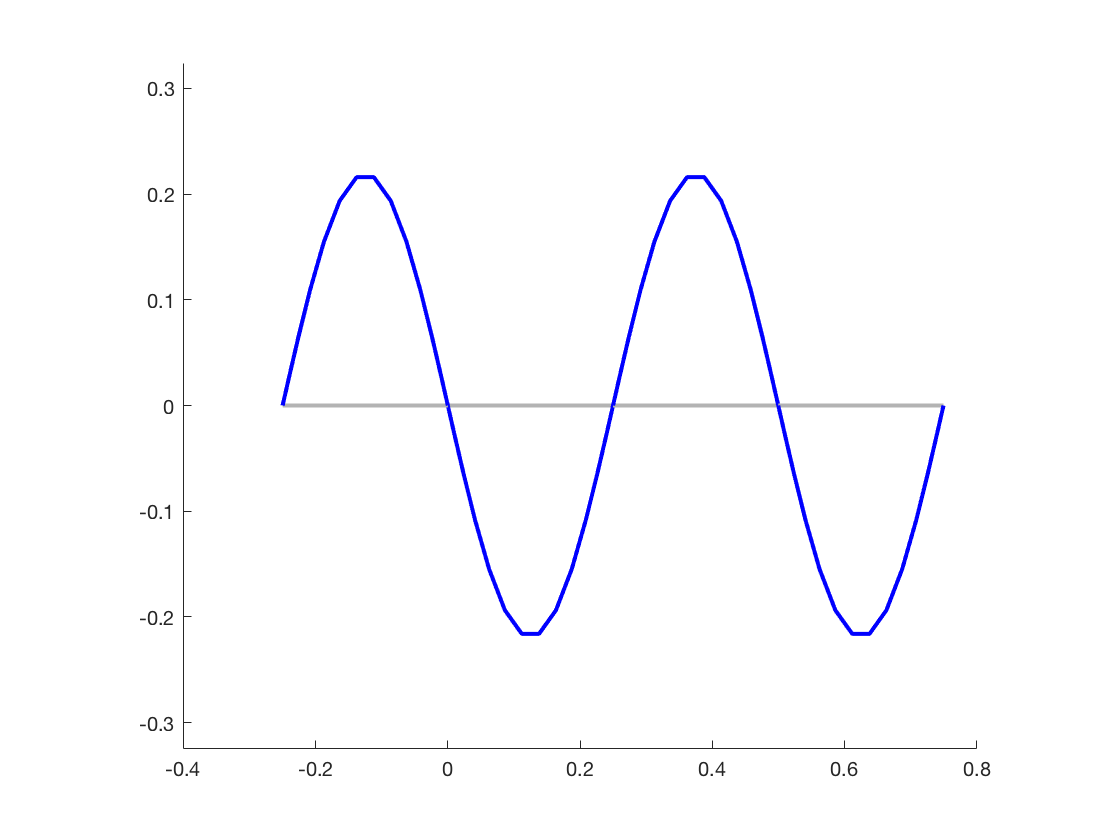}
    \caption{(Top Left) $u_{1,4}(x)$, $\lambda_1=157.9132$. (Top Right) $u_{2,4}(x)$, $\lambda_2=157.9134$.  (Bottom Left) $u_{3,4}(x)$, $\lambda_3=157.9136$. (Bottom Right)  $u_{4,4}(x)$, $\lambda_4=157.9137$.}
    \label{fig3}
\end{figure}

Based on the eigenvalues above and the $M_{\rm thm}$ and $M_{\rm norm}$ matrices for $n=4$, we obtain at $\sigma = 10^7$ the values:
\begin{center}
\begin{tabular}{ |p{3cm}|p{3cm}| }
 \hline
 diff$_{\rm val}$ & {diff}$_{\rm vec}$ \\
 \hline
$ 2.7314 * 10^{-6}$ & $2.7492 * 10^{-7}$ \\
 \hline
\end{tabular}
\end{center}
with the eigenfunctions at $\sigma = 10^7$ pictured in Figure \ref{fig3}. \\

\paragraph{($n = 5$)}

The eigenvalues for five subintervals with increasing $\sigma$ are:

\begin{center}
\begin{tabular}{|p{1cm}|p{1.5cm}|p{1.5cm}|p{1.5cm}|p{1.5cm}|p{1.5cm}| }
 \hline
  & $\sigma = 0$ & $\sigma = 10$ & $\sigma = 10^3$ & $\sigma = 10^5$ & $\sigma = 10^7$ \\
 \hline
 $\lambda_1$ & 9.8696 & 52.3588 & 238.0524 & 246.6509 & 246.7392 \\
 $\lambda_2$ & 39.4784 & 81.3093 & 240.4072 & 246.6755 & 246.7395 \\
 $\lambda_3$ & 88.8264 & 129.0663 & 243.3661 & 246.7060 & 246.7398 \\
 $\lambda_4$ & 157.9137 & 193.3869 & 245.8004 & 246.7307 & 246.7400 \\
 $\lambda_5$ & 246.7401 & 246.7401 & 246.7401 & 246.7401 & 246.7401 \\
 \hline
\end{tabular}
\end{center}

\begin{figure}[ht]
  \centering
    \includegraphics[width=.3\textwidth]{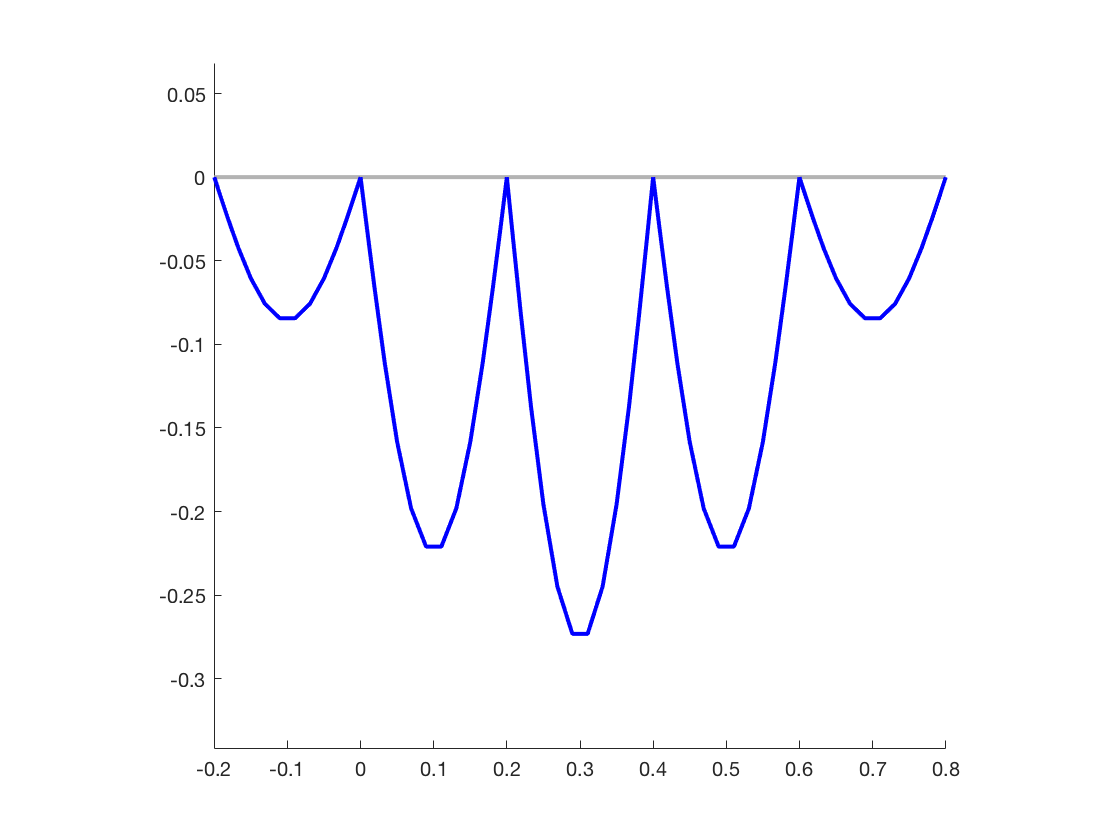}
    \includegraphics[width=.3\textwidth]{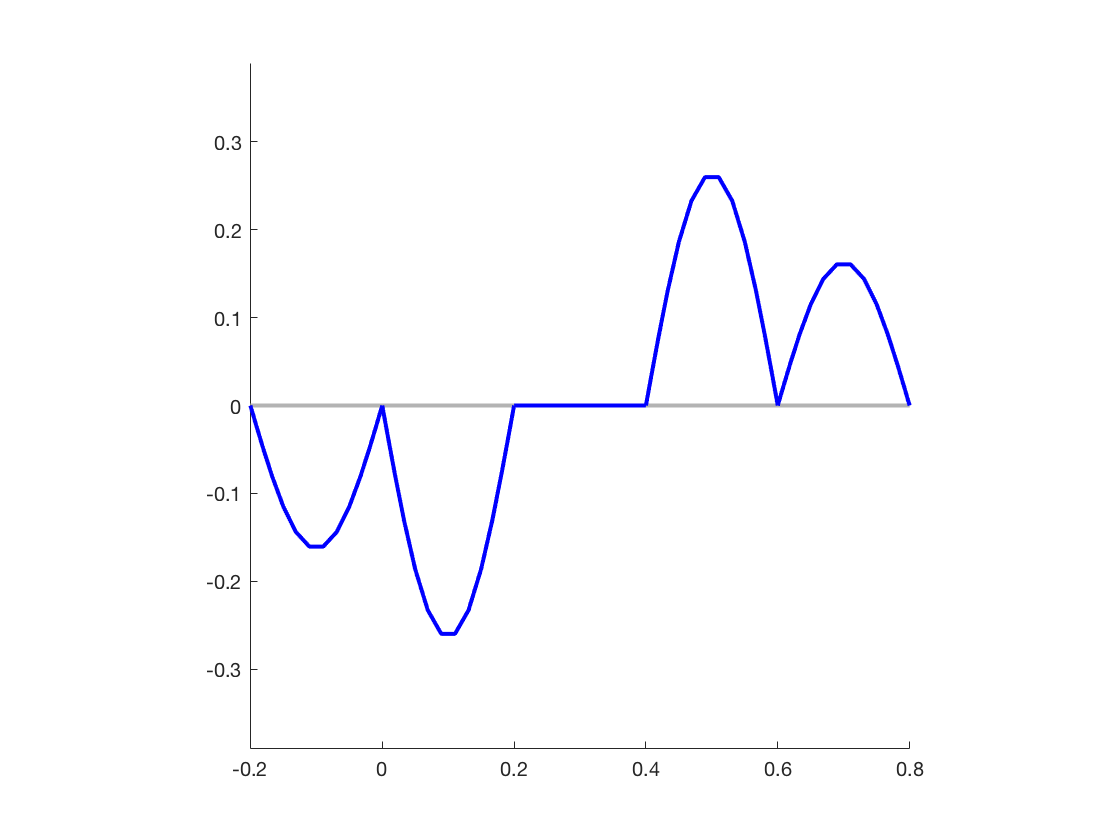}
    \includegraphics[width=.3\textwidth]{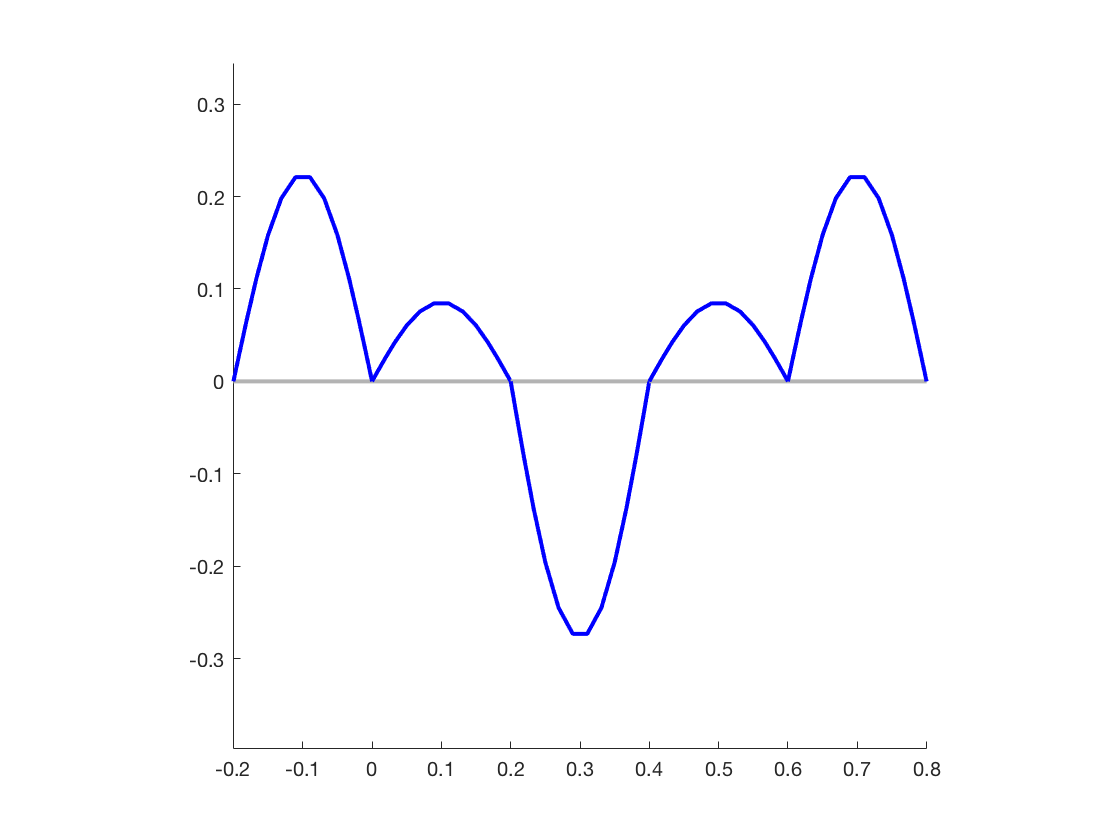} \\
    \includegraphics[width=.3\textwidth]{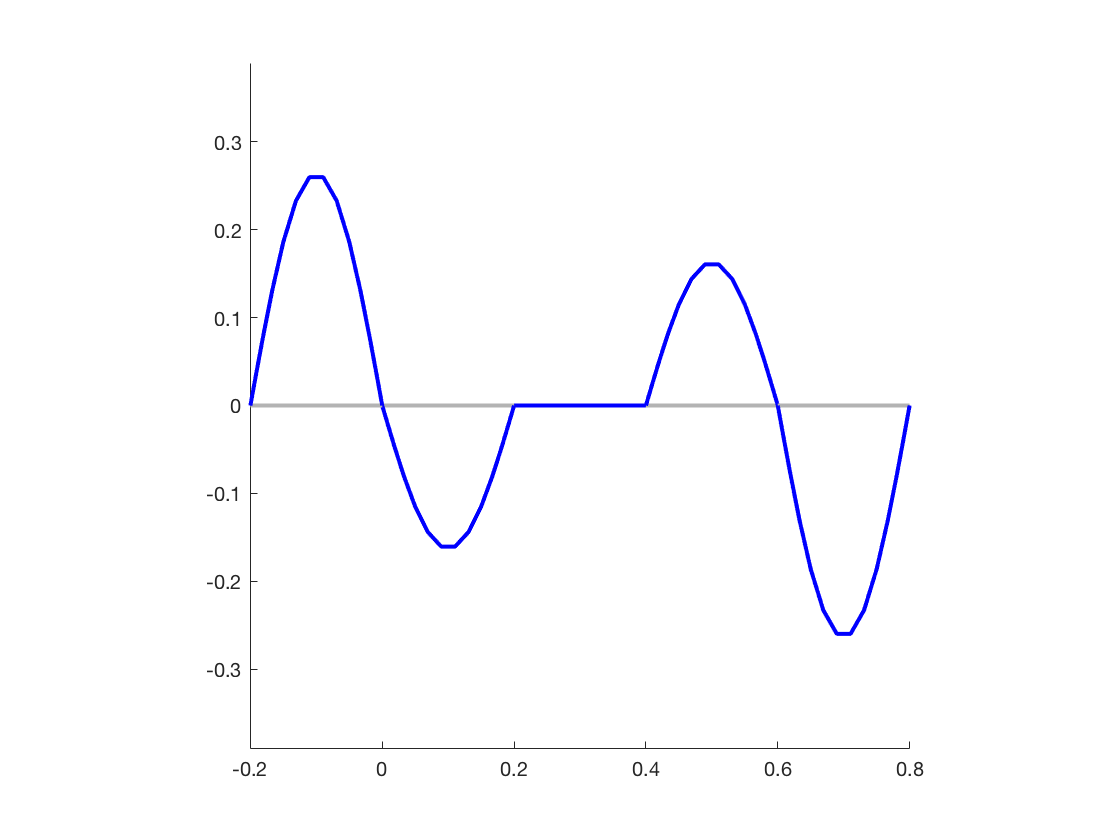}
    \includegraphics[width=.3\textwidth]{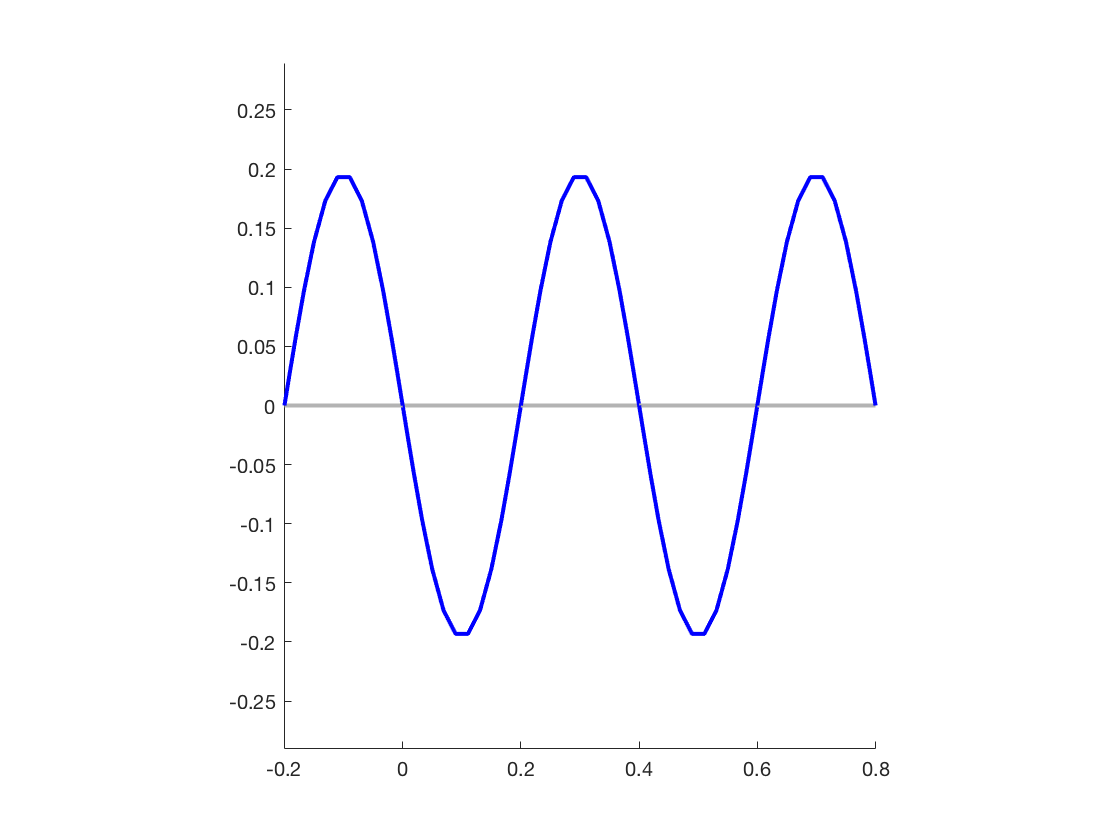}
    \caption{(Top Left) $u_{1,5}(x)$, $\lambda_1=246.7392$. (Top Center) $u_{2,5}(x)$, $\lambda_2=246.7395$. (Top Right) $u_{3,5}(x)$, $\lambda_3=246.7398$. (Bottom Left) $u_{4,5}(x)$, $\lambda_4=246.7400$. (Bottom Right) $u_{5,5}(x)$, $\lambda_5=246.7401$.}
  \label{fig4}
\end{figure}

The corresponding values for diff$_{\rm vec}$ and diff$_{\rm val}$ are given by
\begin{center}
\begin{tabular}{ |p{3cm}|p{3cm}| }
 \hline
 diff$_{\rm val}$ & {diff}$_{\rm vec}$ \\
 \hline
$ 3.6180 * 10^{-6}$ &  $1.5143 * 10^{-6}$ \\
 \hline
 \end{tabular}
\end{center}
and the graphs of the first five eigenfunctions at $\sigma = 10^7$ are pictured in Figure \ref{fig4}. \\

\paragraph{($n = 6$)}

\begin{figure}[ht]
  \centering
    \includegraphics[width=.3\textwidth]{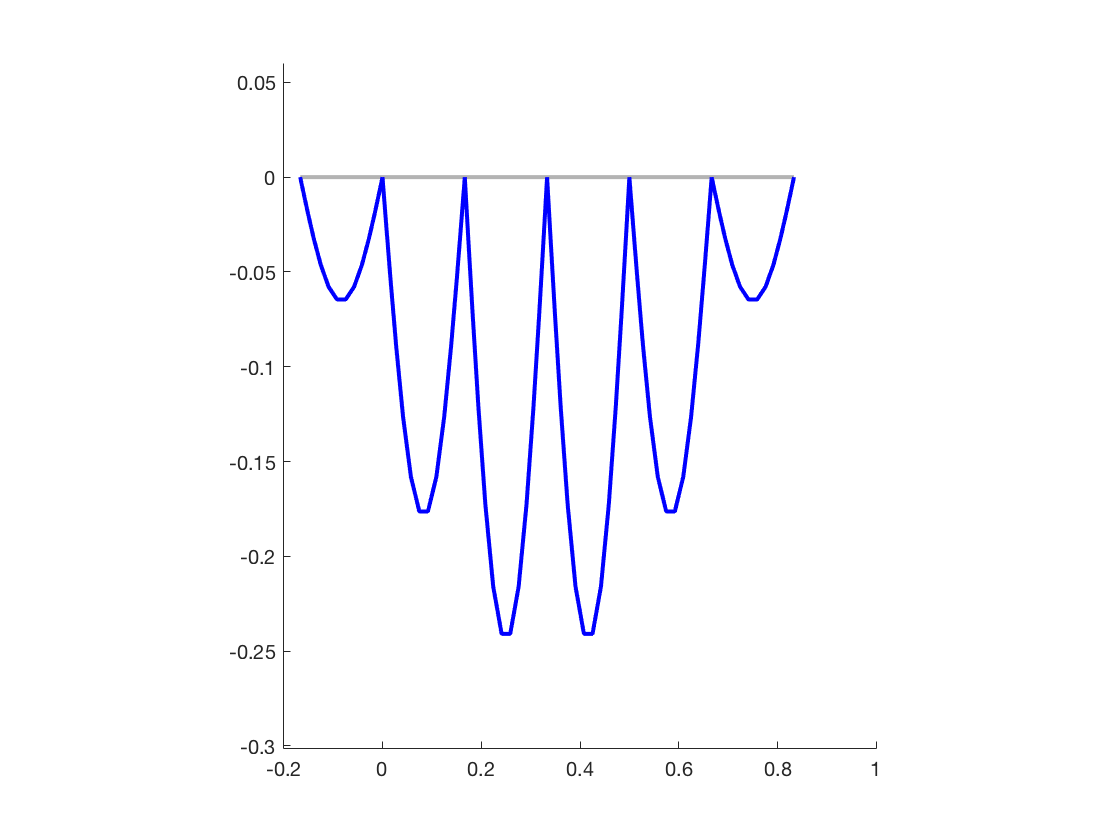}
    \includegraphics[width=.3\textwidth]{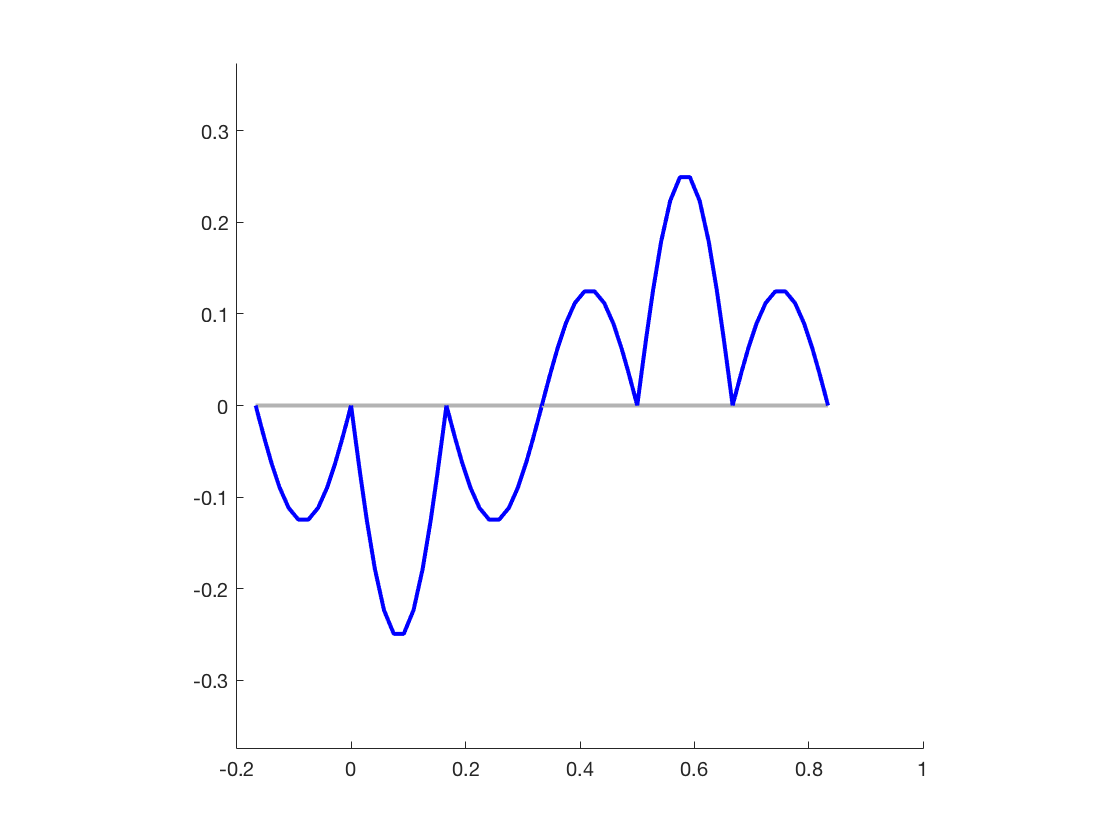}
    \includegraphics[width=.3\textwidth]{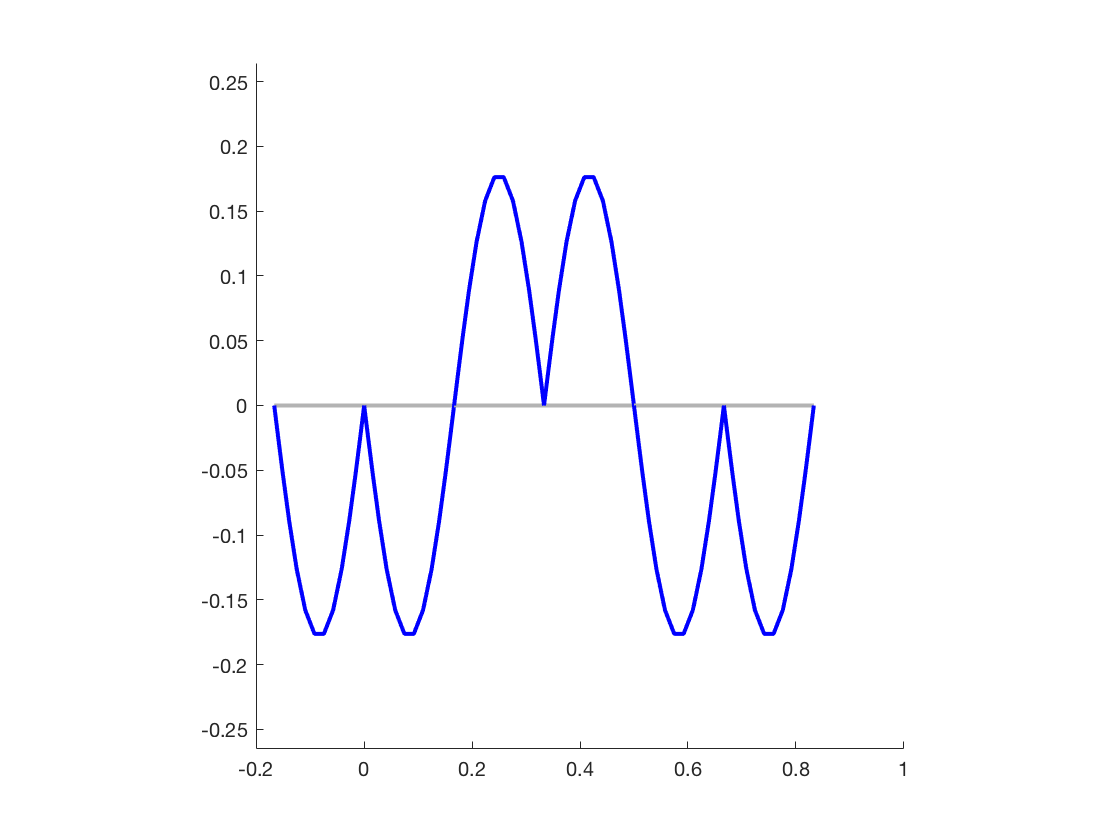} \\
    \includegraphics[width=.3\textwidth]{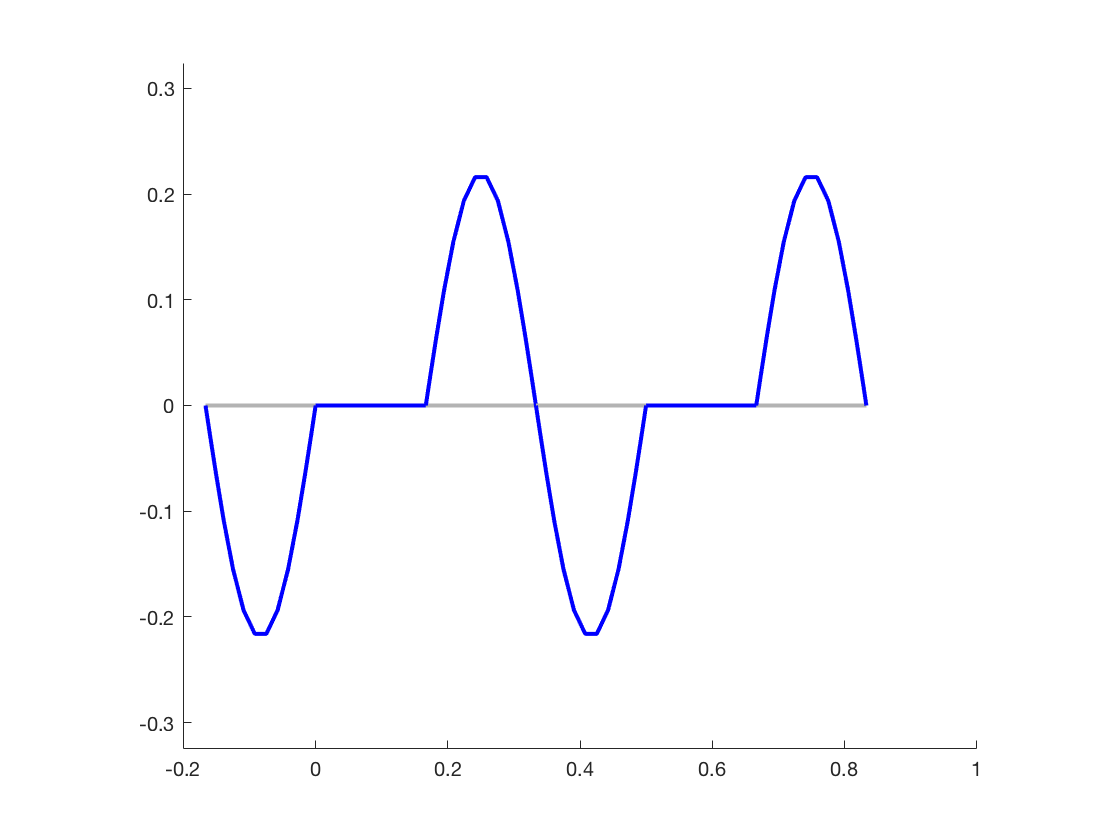}
    \includegraphics[width=.3\textwidth]{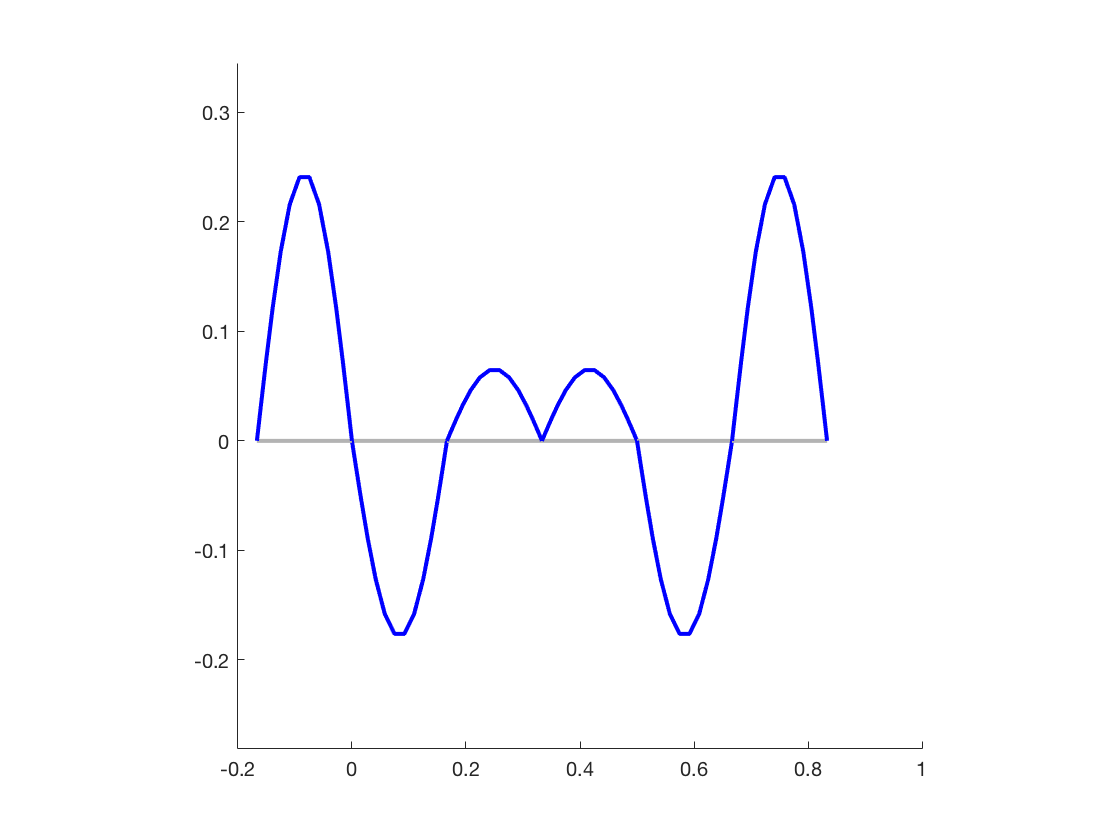}
    \includegraphics[width=.3\textwidth]{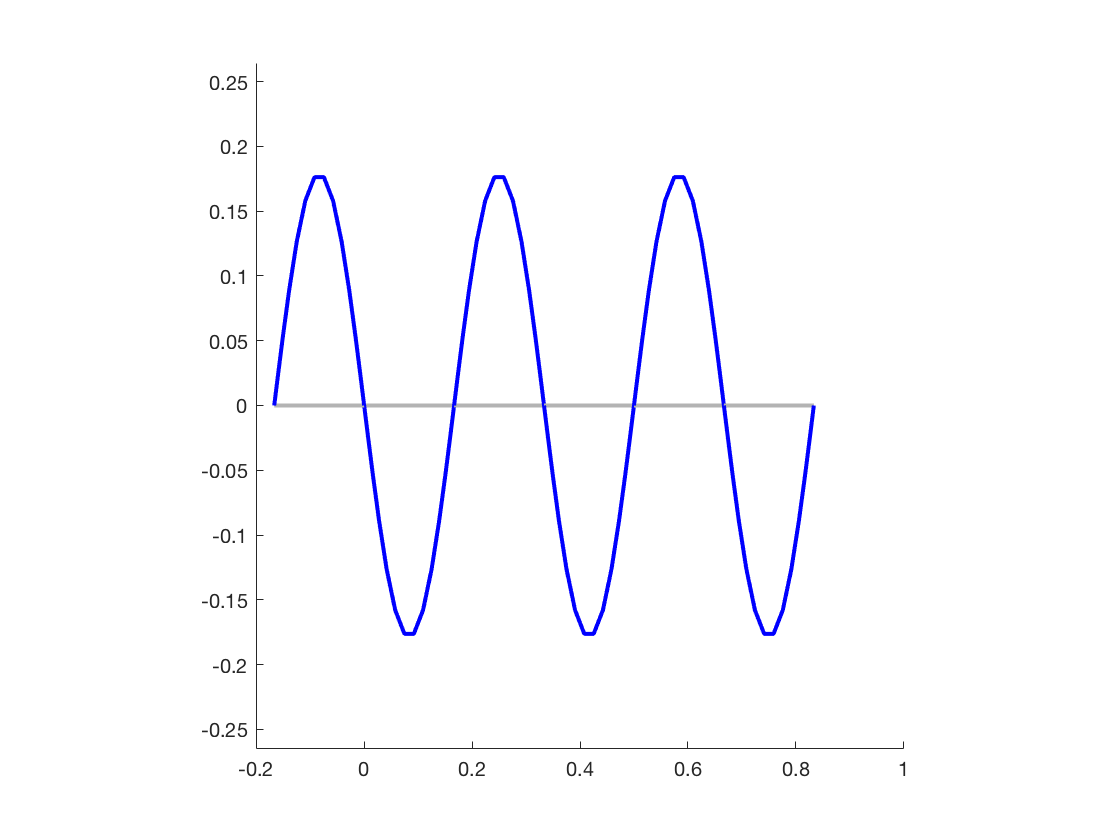}
    \caption{(Top Left) $u_{1,6}(x)$, $\lambda_1=355.3042$. (Top Center) $u_{2,6}(x)$, $\lambda_2=355.3045$. (Top Right) $u_{3,6}(x)$, $\lambda_3=355.3049$. (Bottom Left) $u_{4,6}(x)$, $\lambda_4=355.3053$. (Bottom Center) $u_{5,6}(x)$, $\lambda_5=355.3056$. (Bottom Right) $u_{6,6}(x)$, $\lambda_6=355.3058$.}
  \label{fig5}
\end{figure}

For six subintervals, we obtain the eigenfunctions displayed in Figure \ref{fig5}, with the following and diff$_{\rm val}$ and diff$_{\rm vec}$ values when $\sigma = 10^7$:
\begin{center}
\begin{tabular}{ |p{3cm}|p{3cm}| }
 \hline
 diff$_{\rm val}$ & {diff}$_{\rm vec}$ \\
 \hline
 $4.4784 * 10^{-6}$ & $1.6547 *10^{-6}$ \\
 \hline
\end{tabular}
\end{center}

\appendix

\section{Comparison to Chebfun Method}
\label{app:chebfun}

An additional method for finding the eigenvalues and eigenfunctions as $\sigma$ approaches $\infty$ was implemented using the MATLAB package Chebfun. The Chebfun system allows one to solve differential equations in one dimension using simple operations; more information on Chebfun can be found in \cite{chebfun}.   

This method was implemented for two and three subintervals on $[0,1]$. The general strategy is to define an operator as a vector of $n$ second-order differential operators, one on each of the $n$ subintervals, one vector for the left-side boundary conditions, and one vector for the right-side boundary conditions. The eigenvalues of the operator are then calculated using the \mcode{eigs} command. 

\subsection{Results}
\paragraph{($n = 2$)}
For two subintervals, we write the eigenfunction $u(x)$ as 
\[
    u(x) = \begin{cases} 
    w(x),   &  x \in [0, 1/2]\\
      v(1-x),  & x \in [1/2, 1].
      \end{cases}
\]
 Here we halve the interval and define the eigenfunction $u(x)$ via the two functions $w(x)$ and $v(x)$, in order to apply the boundary conditions at the interior node in addition to the Dirchlet conditions at the endpoints. Specifically, to use Chebfun's $\mcode{L.lbc}$ and $\mcode{L.rbc}$, we must reflect $v(x)$ across the vertical axis in order for $x=0$ to correspond to the left boundary and $x=1/2$ to the right boundary of both $w(x)$ and $v(x)$.
 
As expected, when $\sigma = 0$, the $\mcode{Chebfun}$ version of the $\mcode{eigs}$ operation returns the eigenvalues of \eqref{diffeq}, $\lambda_1 = 9.8696$ and $\lambda_2 = 39.4784$. As $\sigma$ is increased, the first eigenvalue converges to the value of the constant second eigenvalue. For example, for increasing values of $\sigma$, $\mcode{eigenvalues}$ returns the values in the following table:

\begin{center}
\begin{tabular}{ |p{1.5cm}|p{1.5cm}|p{1.75cm}|p{1.75cm}|  }
 \hline
  & $\sigma = 500$ & $\sigma = 1000$ & $\sigma = 5000$ \\
 \hline
 $\lambda_1$ & 38.8544    & 39.1645&   39.4153\\
 $\lambda_2$ &  39.4784  & 39.4784   & 39.4784\\
 \hline
\end{tabular}
\end{center}
The vector corresponding to the $\sigma=\infty$ limit of the  first eigenfunction (see Corollary \ref{cor:infty}) consists of the coefficients of $\sin(2\pi x)$ on each of the two subintervals. This eigenvector can be determined by finding the value of the function at the midpoints of each subinterval.

The signed maxima on each subinterval can then be compared to the eigenvector from Corollary~\ref{cor:infty} as $\sigma$ becomes large. When $\sigma = 5000$, the computation returns the vector $(1.4136,1.4136)$ for the corresponding maxima.  The eigenfunction $\sin(2 \pi x)$ is positive on $(0, 1/2)$ and negative on $(1/2, 1)$. Thus, normalizing the eigenvector, we can see that it is close to $(1,-1)$,  
which is consistent with Corollary \ref{cor:infty} and also with the solution constructed above in Section \ref{sec:num}. \\

\paragraph{($n = 3$)} For three subintervals, we follow a similar process as for $n = 2$, where $u(x)$ is defined as 
\begin{equation*}
     u(x) = \begin{cases}
      w(x),   &  x \in [0, 1/3]\\
      v(2/3-x),  & x \in [1/3, 2/3]\\
      z(x-2/3), & x \in [2/3, 1].
    \end{cases}
\end{equation*}
One can then define the operator and boundary conditions appropriately on the interval $[0,1/3]$, taking care with regards to sign orientation of the derivative on each interval.  Note that the interval of interest is now $[0, 1/3]$; in general, for $n$ subintervals, the operator is applied to the interval $[0, 1/n]$.

As in the $n=2$ case, the first three eigenvalues at $\sigma = 0$ are within machine precision of $\pi^2$, $4 \pi^2$ and $9 \pi^2$.  We then obtain the following for the first three eigenvalues of $L$ with increasing $\sigma$:

\begin{center}
\begin{tabular}{ |p{1.5cm}|p{1.5cm}|p{1.75cm}|p{1.75cm}|  }
 \hline
  & $\sigma = 500$ & $\sigma = 1000$ & $\sigma = 5000$ \\
 \hline
 $\lambda_1$ & 85.7146 & 87.2491 &   88.5075\\
 $\lambda_2$ &  87.7703 & 88.2959  & 88.7199\\
 $\lambda_3$ & 88.8264 & 88.8264 & 88.8264\\
 \hline
\end{tabular}
\end{center}
To find the vectors of the amplitudes of the first and second eigenfunctions, we construct the eigenfunctions as above and return the amplitudes on each interval.

We obtain the vector $(1,-2,1)$ for the first eigenfunction and the vector $(1, 0 ,-1)$ for the second eigenfunction, which are consistent with the vectors in Corollary \ref{cor:infty}, and also matches nicely with our computations in Section \ref{sec:num}.

\bibliographystyle{alpha}
\bibliography{References}

\begin{thebibliography}{HHOT10}

\bibitem[BCM19]{nodal}
Gregory Berkolaiko, Graham Cox, and Jeremy~L Marzuola.
\newblock Nodal deficiency, spectral flow, and the {D}irichlet-to-{N}eumann
  map.
\newblock {\em Letters in Mathematical Physics}, 109(7):1611--1623, 2019.

\bibitem[Ber17]{graph}
Gregory Berkolaiko.
\newblock An elementary introduction to quantum graphs.
\newblock {\em Geometric and computational spectral theory}, 700:41--72, 2017.

\bibitem[BH16]{berard2016courant}
Pierre B{\'e}rard and Bernard Helffer.
\newblock Courant-sharp eigenvalues for the equilateral torus, and for the
  equilateral triangle.
\newblock {\em Letters in Mathematical Physics}, 106(12):1729--1789, 2016.

\bibitem[BHK20]{berard2020courant}
Pierre B{\'e}rard, Bernard Helffer, and Rola Kiwan.
\newblock Courant-sharp property for {D}irichlet eigenfunctions on the
  {M}\"obius strip.
\newblock {\em arXiv preprint arXiv:2005.01175}, 2020.

\bibitem[BKS12]{BKS12}
Gregory Berkolaiko, Peter Kuchment, and Uzy Smilansky.
\newblock Critical partitions and nodal deficiency of billiard eigenfunctions.
\newblock {\em Geom. Funct. Anal.}, 22(6):1517--1540, 2012.

\bibitem[BNH17]{bonnaillie2015nodal}
Virginie Bonnaillie-No\"{e}l and Bernard Helffer.
\newblock Nodal and spectral minimal partitions---the state of the art in 2016.
\newblock In {\em Shape optimization and spectral theory}, pages 353--397. De
  Gruyter Open, Warsaw, 2017.

\bibitem[CJM17]{CJM_nodal}
Graham Cox, Christoper K. R.~T. Jones, and Jeremy~L. Marzuola.
\newblock Manifold decompositions and indices of {S}chr\"{o}dinger operators.
\newblock {\em Indiana Univ. Math. J.}, 66(5):1573--1602, 2017.

\bibitem[DBT08]{chebfun}
Tobin~A Driscoll, Folkmar Bornemann, and Lloyd~N Trefethen.
\newblock The chebop system for automatic solution of differential equations.
\newblock {\em BIT Numerical Mathematics}, 48(4):701--723, 2008.

\bibitem[DH16]{hale}
Tobin~A Driscoll and Nicholas Hale.
\newblock Rectangular spectral collocation.
\newblock {\em IMA Journal of Numerical Analysis}, 36(1):108--132, 2016.

\bibitem[Goo19]{roy}
Roy~H. Goodman.
\newblock N{LS} bifurcations on the bowtie combinatorial graph and the dumbbell
  metric graph.
\newblock {\em Discrete Contin. Dyn. Syst.}, 39(4):2203--2232, 2019.

\bibitem[HHOT09]{helffer2009nodal}
Bernard Helffer, Thomas Hoffmann-Ostenhof, and Susanna Terracini.
\newblock Nodal domains and spectral minimal partitions.
\newblock In {\em Annales de l'IHP Analyse non lin{\'e}aire}, volume~26, pages
  101--138, 2009.

\bibitem[HHOT10]{helffer2010spectral}
Bernard Helffer, Thomas Hoffmann-Ostenhof, and Susanna Terracini.
\newblock On spectral minimal partitions: the case of the sphere.
\newblock In {\em Around the Research of Vladimir Maz'ya III}, pages 153--178.
  Springer, 2010.

\bibitem[HS16]{helffer2016nodal}
Bernard Helffer and Mikael Sundqvist.
\newblock On nodal domains in {E}uclidean balls.
\newblock {\em Proceedings of the American Mathematical Society},
  144(11):4777--4791, 2016.

\bibitem[L{\'e}n15]{lena2015courant}
Corentin L{\'e}na.
\newblock Courant-sharp eigenvalues of a two-dimensional torus.
\newblock {\em Comptes Rendus Mathematique}, 353(6):535--539, 2015.

\bibitem[Ple56]{P56}
{\AA}ke Pleijel.
\newblock Remarks on {C}ourant's nodal line theorem.
\newblock {\em Comm. Pure Appl. Math.}, 9:543--550, 1956.

\bibitem[XH16]{xu}
Kuan Xu and Nicholas Hale.
\newblock Explicit construction of rectangular differentiation matrices.
\newblock {\em IMA Journal of Numerical Analysis}, 36(2):618--632, 2016.

\end{thebibliography}

\end{document}